\documentclass[11pt,fleqn]{article}%
\usepackage{amsfonts}
\usepackage{amsmath}
\usepackage{amssymb}
\usepackage{graphicx}
\usepackage{pdfsync}%
\providecommand{\U}[1]{\protect\rule{.1in}{.1in}}
\newtheorem{theorem}{Theorem}
\newtheorem{assump}{Assumption}

\newtheorem{lemma}[theorem]{Lemma}

\newtheorem{proposition}[theorem]{Proposition}
\newtheorem{remark}[theorem]{Remark}

\newenvironment{proof}[1][Proof]{\noindent\textbf{#1.} }{\ \rule{0.5em}{0.5em}}

\numberwithin{equation}{section}
\numberwithin{theorem}{section}
\begin{document}

\title{Adaptive estimation of spectral densities via wavelet thresholding and
information projection}
\author{J\'{e}r\'{e}mie Bigot$^{1}$, Rolando J. Biscay Lirio$^{3}$, Jean-Michel
Loubes$^{1}$ \& Lilian Mu\~{n}iz Alvarez$^{1,2}$\\IMT, Universit\'{e} Paul Sabatier, Toulouse, France$^{1}$\\Facultad de Matem\'{a}tica y Computaci\'{o}n, Universidad de La Habana,
Cuba$^{2}$\\CIMFAV-DEUV, Facultad de Ciencias, Universidad de Valparaiso, Chile$^{3}$}
\maketitle

\begin{abstract}
In this paper, we study the problem of adaptive estimation of the spectral
density of a stationary Gaussian process. For this purpose, we consider a
wavelet-based method which combines the ideas of wavelet approximation and
estimation by information projection in order to warrants that the solution is
a non-negative function. The spectral density of the process is estimated by
projecting the wavelet thresholding expansion of the periodogram onto a family
of exponential functions. This ensures that the spectral density estimator is
a strictly positive function. The theoretical behavior of the estimator is
established in terms of rate of convergence of the Kullback-Leibler
discrepancy over Besov classes. We also show the excellent practical
performance of the estimator in some numerical experiments.

\end{abstract}

\thispagestyle{empty}



\section{Introduction}

The estimation of spectral densities is a fundamental problem in inference for
stationary stochastic processes. Many applications in several fields such as
weather forecast and financial series are deeply related to this issue, see
for instance Priestley \cite{Pr}. It is known that the estimation of the
covariance function of a stationary process is strongly related to the
estimation of the corresponding spectral density. By Bochner's theorem the
covariance function is non-negative definite if and only if the corresponding
spectral density is a non-negative function. Hence in order to preserve the
property of non-negative definiteness of a covariance function, the estimation
of the corresponding spectral density must be a non-negative function. The
purpose of this work is to provide a non-negative estimator of the spectral density.

Inference in the spectral domain uses the periodogram of the data, providing
an inconsistent estimator which must be smoothed in order to achieve
consistency. For highly regular spectral densities, linear smoothing
techniques such as kernel smoothing are appropriate (see Brillinger
\cite{Br}). However, these methods are not able to achieve the
optimal mean-square rate of convergence for spectra whose smoothness is
distributed inhomogeneously over the domain of interest. For this nonlinear
methods are needed. One nonlinear method for adaptive spectral density
estimation of a stationary Gaussian sequence was proposed by Comte
\cite{Comte}. It is based on model selection techniques. Others nonlinear
smoothing procedures are the wavelet thresholding methods, first proposed by
Donoho and Johnstone \cite{D-J}. In this context, different thresholding rules
have been proposed by Neumann \cite{N} and Fryzlewicz, Nason and von Sachs
\cite{F} to name but a few.

Neumann's approach \cite{N} consists in pre-estimating the variance of the
periodogram via kernel smoothing, so that it can be supplied to the wavelet
estimation procedure. Kernel pre-estimation may not be appropriate in cases
where the underlying spectral density is of low regularity. One way to avoid
this problem is proposed in Fryzlewicz, Nason and von Sachs \cite{F}, where
the empirical wavelet coefficient thresholds are built as appropriate local
weighted $l_{1}$ norms of the periodogram. These methods do not produce
non-negative spectral density estimators, therefore the corresponding
estimators of the covariance function is not non-negative definite.

To overcome the drawbacks of previous estimators, in this paper we propose a
new wavelet-based method for the estimation of the spectral density of a
Gaussian process. As a solution to ensure non-negativeness of the spectral
density estimator, our method combines the ideas of wavelet thresholding and
estimation by information projection. We estimate the spectral density by a
projection of the nonlinear wavelet approximation of the periodogram onto a
family of exponential functions. Therefore, the estimator is non-negative by
construction. This technique was studied by Barron and Sheu \cite{B-S} for the
approximation of density functions by sequences of exponential families, by
Loubes and Yan \cite{L-Y} for penalized maximum likelihood estimation with
$l_{1}$ penalty, by Antoniadis and Bigot \cite{A-B} for the study of Poisson
inverse problems, and by Bigot and Van Bellegem \cite{B-V} for log-density deconvolution.

The theoretical optimality of the estimators for the spectral density of a
stationary process is generally studied using risk bounds in $L_{2}$-norm.
This is the case in the papers of Neumann \cite{N}, Comte \cite{Comte} and
Fryzlewicz, Nason and von Sachs \cite{F} mentioned before. In this work, the
behavior of the proposed estimator is established in terms of the rate of
convergence of the Kullback-Leibler discrepancy over Besov classes, which is
maybe a more natural loss function for the estimation of a spectral density
function than the $L_{2}$-norm. Moreover, the thresholding rules that we use
to derive adaptive estimators differ from previous approaches based on wavelet
decomposition and are quite simple to compute. Finally, we compare the
performance of our estimator with other estimators on some simulations.

The paper is organized as follows. Section 2 presents the statistical
framework under which we work. We define the model, the wavelet-based
exponential family and the linear and nonlinear wavelet estimators by
information projection. We also recall the definition of the Kullback-Leibler
divergence and some results on Besov spaces. The rate of convergence of the
proposed estimators are stated in Section 3. Some numerical experiments are
described in Section 4. Technical lemmas and proofs of the main theorems are
gathered in the Appendix.

Throughout this paper $C$ denotes a constant that may vary from line to line.
The notation $C(.)$ specifies the dependency of $C$ on some quantities.

\section{Statistical framework}

\subsection{The model}

We aim at providing a nonparametric adaptive estimation of the spectral
density which satisfies the property of being non-negative in order to
guarantee that the covariance estimator is a non-negative definite function.
We consider the sequence $\left(  X_{t}\right)  _{t\in\mathbb{N}}$ that
satisfies the following assumptions:

\begin{assump}
\label{ass1} The sequence $\left(  X_{1},...X_{n}\right)  $ is an $n$-sample
drawn from a stationary sequence of Gaussian random variables.
\end{assump}

Let $\rho$ be the covariance function of the process, i.e. $\rho\left(
h\right)  =cov\left(  X_{t},X_{t+h}\right)  $ with $h\in\mathbb{Z}$. The
spectral density $f$ is defined as:%
\[
f\left(  \omega\right)  =\frac{1}{2\pi}\sum\limits_{h\in\mathbb{Z}}\rho\left(
h\right)  e^{i2\pi\omega h},\;\omega\in\left[  0,1\right]  .
\]
We need the following standard assumption on $\rho$:

\begin{assump}
\label{ass2} The covariance function $\rho$ is non-negative definite, such
that there exists two constants $0<C_{1},C_{2}<+\infty$ such that
$\sum\limits_{h\in\mathbb{Z}}\left\vert \rho\left(  h\right)  \right\vert
=C_{1}$ and $\sum\limits_{h\in\mathbb{Z}}\left\vert h\rho^{2}\left(  h\right)
\right\vert =C_{2}$.
\end{assump}

Assumption 2 implies in particular that the spectral density $f$ is bounded by
the constant $C_{1}$. As a consequence, it is also square integrable. As in
Comte \cite{Comte}, the data consist on a number of observations
$X_{1},...,X_{n}$ at regularly spaced points. We want to obtain a positive
estimator for the spectral density function $f$ without parametric assumptions
on the basis of these observations. For this, we combine the ideas of wavelet
thresholding and estimation by information projection.

\subsection{Estimation by information projection}

To ensure nonnegativity of the estimator, we will look for approximations over
an exponential family. For this, we construct a sieve of exponential functions
defined in a wavelet basis.\vskip.1in Let $\phi\left(  \omega\right)  $ and
$\psi\left(  \omega\right)  $, respectively, be the scaling and the wavelet
functions generated by an orthonormal multiresolution decomposition of
$L_{2}\left(  \left[  0,1\right]  \right)  $, see Mallat \cite{Ma} for a
detailed exposition on wavelet analysis. Throughout the paper, the functions
$\phi$ and $\psi$ are supposed to be compactly supported and such that
$\left\Vert \phi\right\Vert _{\infty}<+\infty$, $\left\Vert \psi\right\Vert
_{\infty}<+\infty$. Then, for any integer $j_{0}\geq0$, any function $g\in
L_{2}\left(  \left[  0,1\right]  \right)  $ has the following representation:%
\[
g\left(  \omega\right)  =\sum\limits_{k=0}^{2^{j_{0}}-1}\left\langle
g,\phi_{j_{0},k}\right\rangle \phi_{j_{0},k}\left(  \omega\right)
+\sum\limits_{j=j_{0}}^{+\infty}\sum\limits_{k=0}^{2^{j}-1}\left\langle
g,\psi_{j,k}\right\rangle \psi_{j,k}\left(  \omega\right)  ,
\]
where $\phi_{j_{0},k}\left(  \omega\right)  =2^{\frac{j_{0}}{2}}\phi\left(
2^{j_{0}}\omega-k\right)  $ and $\psi_{j,k}\left(  \omega\right)  =2^{\frac
{j}{2}}\psi\left(  2^{j}\omega-k\right)  $. The main idea of this paper is to
expand the spectral density $f$ onto this wavelet basis and to find an
estimator of this expansion that is then modified to impose the positivity
property. The scaling and wavelet coefficients of the spectral density
function $f$ are denoted by $a_{j_{0},k}=\left\langle f,\phi_{j_{0}%
,k}\right\rangle $ and $b_{j,k}=\left\langle f,\psi_{j,k}\right\rangle $.

To simplify the notations, we write $\left(  \psi_{j,k}\right)  _{j=j_{0}-1}$
for the scaling functions $\left(  \phi_{j,k}\right)  _{j=j_{0}}$. Let
$j_{1}\geq j_{0}$ and define the set
\[
\Lambda_{j_{1}}=\left\{  \left(  j,k\right)  :j_{0}-1\leq j<j_{1},0\leq
k\leq2^{j}-1\right\}  .
\]
Note that $\#\Lambda_{j_{1}}=2^{j_{1}}$, where $\#\Lambda_{j_{1}}$ denotes the
cardinal of $\Lambda_{j_{1}}$. Let $\theta$ denotes a vector in $\mathbb{R}%
^{\#\Lambda_{j_{1}}}$, the wavelet-based exponential family $\mathcal{E}%
_{j_{1}}$ at scale $j_{1}$ is defined as the set of functions:%
\begin{equation}
\mathcal{E}_{j_{1}}=\left\{  f_{j_{1},\theta}\left(  .\right)  =\exp\left(
\sum\limits_{(j,k)\in\Lambda_{j_{1}}}\theta_{j,k}\psi_{j,k}\left(  .\right)
\right)  ,\;\theta=\left(  \theta_{j,k}\right)  _{(j,k)\in\Lambda_{j_{1}}}%
\in\mathbb{R}^{\#\Lambda_{j_{1}}}\right\}  . \label{ExpFamily}%
\end{equation}

It is well known that Besov spaces for periodic functions in $L_{2}([0,1])$
can be characterized in terms of wavelet coefficients (see e.g. Mallat
\cite{Ma}). Assume that $\psi$ has $m$ vanishing moments, and let $0<s<m$
denote the usual smoothness parameter. Then, for a Besov ball $B_{p,q}^{s}(A)$
of radius $A>0$ with $1\leq p,q\leq\infty$, one has that for $s^{\ast
}=s+1/2-1/p\geq0$:
\[
B_{p,q}^{s}\left(  A\right)  :=\left\{  g\in L_{2}([0,1]:\Vert g\Vert
_{s,p,q}:=\left(  \sum\limits_{k=0}^{2^{j_{0}}-1}\left\vert a_{j_{0}%
,k}\right\vert ^{p}\right)  ^{\frac{1}{p}}+\left(  \sum\limits_{j=j_{0}%
}^{\infty}2^{js^{\ast}q}\left(  \sum\limits_{k=0}^{2^{j}-1}\left\vert
b_{j,k}\right\vert ^{p}\right)  ^{\frac{q}{p}}\right)  ^{\frac{1}{q}}\leq
A\right\}  ,
\]
with the respective above sums replaced by maximum if $p=\infty$ or $q=\infty$
and where $a_{j_{0},k}=\left\langle g,\phi_{j_{0},k}\right\rangle $ and
$b_{j,k}=\left\langle g,\psi_{j,k}\right\rangle $.

The condition that $s+1/2-1/p\geq0$ is imposed to ensure that $B_{p,q}^{s}(A)$
is a subspace of $L_{2}([0,1])$, and we shall restrict ourselves to this case
in this paper (although not always stated, it is clear that all our results
hold for $s<m$). 

Let $M>0$ and denote by $F_{p,q}^{s}(M)$ the set of functions such that%
\[
F_{p,q}^{s}(M)=\left\{  f=\exp\left(  g\right)  :\;\Vert g\Vert_{s,p,q}\leq
M\right\}  ,
\]
where $\Vert g\Vert_{s,p,q}$ denotes the norm in the Besov space $B_{p,q}^{s}%
$. Note that assuming that $f\in F_{p,q}^{s}(M)$ implies that $f$ is strictly
positive. The following results hold.

\begin{lemma}
\label{LemaA.4inAB} Suppose that $f\in F_{p,q}^{s}(M)$ with $s>\frac{1}{p}$
and $1\leq p\leq2$. Then, there exists a constant $M_{1}$ such that for all
$f\in F_{p,q}^{s}(M)$, $0<M_{1}^{-1}\leq f\leq M_{1}<+\infty.$
\end{lemma}

Let $V_{j}$ denote the usual multiresolution space at scale $j$ spanned by the
scaling functions $(\phi_{j,k})_{0\leq k\leq2^{j}-1}$, and define
$A_{j}<+\infty$ as the constant such that $\left\Vert \upsilon\right\Vert
_{\infty}\leq A_{j}\left\Vert \upsilon\right\Vert _{L_{2}}$ for all
$\upsilon\in V_{j}$. For $f\in F_{p,q}^{s}(M)$, let $g=\log\left(  f\right)
$. Then for $j\geq j_{0}-1$, define $D_{j}=\left\Vert g-g_{j}\right\Vert
_{L_{2}}$ and $\gamma_{j}=\left\Vert g-g_{j}\right\Vert _{\infty},$ where
$g_{j}=\sum\limits_{k=0}^{2^{j}-1}\theta_{j,k}\psi_{j,k}$, with $\theta
_{j,k}=\left\langle g,\psi_{j,k}\right\rangle $.

The proof of the following lemma immediately follows from the arguments in the
proof of Lemma A.5 in Antoniadis and Bigot \cite{A-B}.

\begin{lemma}
\label{LemmaAjDj} Let $j\in\mathbb{N}$. Then $A_{j}\leq C2^{j/2}$. Suppose
that $f\in F_{p,q}^{s}(M)$ with $1\leq p\leq2$ and $s>\frac{1}{p}$. Then,
uniformly over $F_{p,q}^{s}(M)$, $D_{j}\leq C2^{-j(s+1/2-1/p)}$ and
$\gamma_{j}\leq C2^{-j(s-1/p)}$ where $C$ denotes constants depending only on
$M$, $s$, $p$ and $q$.
\end{lemma}

To assess the quality of the estimators, we will measure the discrepancy
between an estimator $\widehat{f}$ and the true function $f$ in the sense of
relative entropy (Kullback-Leibler divergence) defined by:%
\[
\Delta\left(  f;\widehat{f}\right)  =\int_{0}^{1}\left(  f\log\left(  \frac
{f}{\widehat{f}}\right)  -f+\widehat{f}\right)  d\mu,
\]
where $\mu$ denotes the Lebesgue measure on $[0,1]$. It can be shown that
$\Delta\left(  f;\widehat{f}\right)  $ is non-negative and equals zero if and
only if $\widehat{f}=f$.

We will enforce our estimator of the spectral density to belong to the family
$\mathcal{E}_{j_{1}}$ of exponential functions, which are positive by
definition. For this we will consider a notion of projection using information projection.

The estimation of density function based on information projection has been
introduced by Barron and Sheu \cite{B-S}. To apply this method in our context,
we recall for completeness a set of results that are useful to prove the
existence of our estimators. The proofs of the following lemmas immediately
follow from results in Barron and Sheu \cite{B-S} and Antoniadis and Bigot
\cite{A-B}.

\begin{lemma}
\label{Lema3.1enAB} Let $\beta\in\mathbb{R}^{\#\Lambda_{j_{1}}}$. Assume that
there exists some $\theta\left(  \beta\right)  \in\mathbb{R}^{\#\Lambda
_{j_{1}}}$ such that, for all $\left(  j,k\right)  \in$ $\Lambda_{j_{1}}$,
$\theta\left(  \beta\right)  $ is a solution of
\[
\left\langle f_{j,\theta\left(  \beta\right)  },\psi_{j,k}\right\rangle
=\beta_{j,k}.
\]
Then for any function $f$ such that $\left\langle f,\psi_{j,k}\right\rangle
=\beta_{j,k}$ for all $\left(  j,k\right)  \in$ $\Lambda_{j_{1}}$, and for all
$\theta\in\mathbb{R}^{\#\Lambda_{j_{1}}}$, the following Pythagorian-like
identity holds:%
\begin{equation}
\Delta\left(  f;f_{j,\theta}\right)  =\Delta\left(  f;f_{j,\theta\left(
\beta\right)  }\right)  +\Delta\left(  f_{j,\theta\left(  \beta\right)
};f_{j,\theta}\right)  . \label{IdPyth}%
\end{equation}

\end{lemma}

The next lemma is a key result which gives sufficient conditions for the
existence of the vector $\theta\left(  \beta\right)  $ as defined in Lemma
\ref{Lema3.1enAB}. This lemma also relates distances between the functions in
the exponential family to distances between the corresponding wavelet
coefficients. Its proof relies upon a series of lemmas on bounds within
exponential families for the Kullback-Leibler divergence and can be found in
Barron and Sheu \cite{B-S} and Antoniadis and Bigot \cite{A-B}.

\begin{lemma}
\label{Lema3.2enAB} Let $\theta_{0}\in\mathbb{R}^{\#\Lambda_{j_{1}}}$,
$\beta_{0}=\left(  \beta_{0,\left(  j,k\right)  }\right)  _{\left(
j,k\right)  \in\Lambda_{j_{1}}}\in\mathbb{R}^{\#\Lambda_{j_{1}}}$ such that
$\beta_{0,\left(  j,k\right)  }=$ $\left\langle f_{j,\theta_{0}},\psi
_{j,k}\right\rangle $ for all $\left(  j,k\right)  \in\Lambda_{j_{1}}$, and
$\widetilde{\beta}\in\mathbb{R}^{\#\Lambda_{j_{1}}}$ a given vector. Let
$b=\exp\left(  \left\Vert \log\left(  f_{j,\theta_{0}}\right)  \right\Vert
_{\infty}\right)  $ and $e=\exp(1)$. If $\left\Vert \widetilde{\beta}%
-\beta_{0}\right\Vert _{2}\leq\frac{1}{2ebA_{j_{1}}}$ then the solution
$\theta\left(  \widetilde{\beta}\right)  $ of
\[
\left\langle f_{j_{1},\theta},\psi_{j,k}\right\rangle =\widetilde{\beta}%
_{j,k}\text{ for all }\left(  j,k\right)  \in\Lambda_{j_{1}}%
\]
exists and satisfies%
\begin{align*}
\left\Vert \theta\left(  \widetilde{\beta}\right)  -\theta_{0}\right\Vert
_{2}  &  \leq2eb\left\Vert \widetilde{\beta}-\beta_{0}\right\Vert _{2}\\
\left\Vert \log\left(  \frac{f_{j_{1},\theta\left(  \beta_{0}\right)  }%
}{f_{j_{1},\theta\left(  \widetilde{\beta}\right)  }}\right)  \right\Vert
_{\infty}  &  \leq2ebA_{j_{1}}\left\Vert \widetilde{\beta}-\beta
_{0}\right\Vert _{2}\\
\Delta\left(  f_{j_{1},\theta\left(  \beta_{0}\right)  };f_{j_{1}%
,\theta\left(  \widetilde{\beta}\right)  }\right)   &  \leq2eb\left\Vert
\widetilde{\beta}-\beta_{0}\right\Vert _{2}^{2},
\end{align*}
where $\Vert\beta\Vert_{2}$ denotes the standard Euclidean norm for $\beta
\in\mathbb{R}^{\#\Lambda_{j_{1}}}$.
\end{lemma}

Following Csisz\'{a}r \cite{Csiszar}, it is possible to define the projection
of a function $f$ onto $\mathcal{E}_{j_{1}}$. If this projection exists, it is
defined as the function $f_{j_{1},\theta_{j_{1}}^{\ast}}$ in the exponential
family $\mathcal{E}_{j_{1}}$ that is the closest to the true function $f$ in
the Kullback-Leibler sense, and is characterized as the unique function in the
family $\mathcal{E}_{j_{1}}$ for which
\[
\left\langle f_{j_{1},\theta_{j_{1}}^{\ast}},\psi_{j,k}\right\rangle
=\left\langle f,\psi_{j,k}\right\rangle :=\beta_{j,k}\mbox{ for all }\left(
j,k\right)  \in\Lambda_{j_{1}}.
\]
Note that the notation $\beta_{j,k}$ is used to denote both the the scaling
coefficients $a_{j_{0},k}$ and the wavelet coefficients $b_{j,k}$.

Let%
\[
I_{n}\left(  \omega\right)  =\frac{1}{2\pi n}\sum\limits_{t=1}^{n}%
\sum\limits_{t^{\prime}=1}^{n}\left(  X_{t}-\overline{X}\right)  \left(
X_{t^{\prime}}-\overline{X}\right)  ^{\ast}e^{i2\pi\omega\left(  t-t^{\prime
}\right)  },
\]
be the classical periodogram, where $\left(  X_{t}-\overline{X}\right)
^{\ast}$ denotes the conjugate transpose of $\left(  X_{t}-\overline
{X}\right)  $ and $\overline{X}=\frac{1}{n}\sum\limits_{t=1}^{n}X_{t}$. The
expansion of $I_{n}\left(  \omega\right)  $ onto the wavelet basis allows to
obtain estimators of $a_{j_{0},k}$ and $b_{j,k}$ given by
\begin{equation}
\widehat{a}_{j_{0},k}=\int\limits_{0}^{1}I_{n}\left(  \omega\right)
\phi_{j_{0},k}\left(  \omega\right)  d\omega\quad\mbox{ and }\quad\widehat
{b}_{j,k}=\int\limits_{0}^{1}I_{n}\left(  \omega\right)  \psi_{j,k}\left(
\omega\right)  d\omega. \label{WCE}%
\end{equation}
It seems therefore natural to estimate the function $f$ by searching for some
$\widehat{\theta}_{n}\in\mathbb{R}^{\#\Lambda_{j_{1}}}$ such that
\begin{equation}
\left\langle f_{j_{1},\widehat{\theta}_{n}},\psi_{j,k}\right\rangle
=\int\limits_{0}^{1}I_{n}\left(  \omega\right)  \psi_{j,k}\left(
\omega\right)  d\omega:=\widehat{\beta}_{j,k}\mbox{ for all }\left(
j,k\right)  \in\Lambda_{j_{1}}, \label{EmpBeta}%
\end{equation}
where $\widehat{\beta}_{j,k}$ denotes both the estimation of the scaling
coefficients $\widehat{a}_{j_{0},k}$ and the wavelet coefficients $\widehat
{b}_{j,k}$. The function $f_{j_{1},\widehat{\theta}_{n}}$ is the spectral
density positive linear estimator.

Similarly, the positive nonlinear estimator with hard thresholding is defined
as the function $f_{j_{1},\widehat{\theta}_{n},\xi}^{HT}$ (with $\widehat
{\theta}_{n}\in\mathbb{R}^{\#\Lambda_{j_{1}}}$) such that
\begin{equation}
\left\langle f_{j_{1},\widehat{\theta}_{n},\xi}^{HT},\psi_{j,k}\right\rangle
=\delta_{\xi}\left(  \widehat{\beta}_{j,k}\right)  \mbox{ for all }\left(
j,k\right)  \in\Lambda_{j_{1}}, \label{ThresBeta}%
\end{equation}
where $\delta_{\xi}$ denotes the hard thresholding rule defined by
\[
\delta_{\xi}\left(  x\right)  =xI\left(  |x|\geq\xi\right)  \mbox{ for }x\in
{\mathbb{R}},
\]
where $\xi>0$ is an appropriate threshold whose choice is discussed later on.

The existence of these estimators is questionable. Thus,
in the next sections, some sufficient conditions are given for the existence
of $f_{j_{1},\widehat{\theta}_{n}}$ and $f_{j_{1},\widehat{\theta}_{n},\xi
}^{HT}$ with probability tending to one as $n\rightarrow+\infty$. Even if an explicit expression for $\widehat{\theta}_{n}$ is not available,  we use a numerical approximation of $\widehat{\theta}_{n}$, obtained via a gradient-descent algorithm with an adaptive step.

In this section we establish the rate of convergence of our estimators in
terms of the Kullback-Leibler discrepancy over Besov classes.

We make the following assumption on the wavelet basis that guarantees that
Assumption \ref{ass2} holds uniformly over $F_{p,q}^{s}(M)$.

\begin{assump}
\label{ass3} Let $M>0$, $1\leq p\leq2$ and $s>1/p$. For $f\in F_{p,q}^{s}(M)$
and $h\in{\mathbb{Z}}$, let $\rho(h)=\int_{0}^{1}f(\omega)e^{-i2\pi\omega
h}d\omega$, $C_{1}(f):=\sum\limits_{h\in\mathbb{Z}}\left\vert \rho\left(
h\right)  \right\vert $ and $C_{2}(f):=\sum\limits_{h\in\mathbb{Z}}\left\vert
h\rho^{2}\left(  h\right)  \right\vert $. Then, the wavelet basis is such that
there exists a constant $M_{\ast}$ such that for all $f\in F_{p,q}^{s}(M)$,
\[
C_{1}(f)\leq M_{\ast}\mbox{ and }C_{2}(f)\leq M_{\ast}.
\]

\end{assump}

\subsection{Linear estimation}

The following theorem is the general result on the linear information
projection estimator of the spectral density function. Note that the choice of
the coarse level resolution level $j_{0}$ is of minor importance, and without
loss of generality we take $j_{0}=0$ for the linear estimator $f_{j_{1}%
,\widehat{\theta}_{n}}$.

\begin{theorem}
\label{theo:lin} Assume that $f\in F_{2,2}^{s}\left(  M\right)  $ with
$s>\frac{1}{2}$ and suppose that Assumptions \ref{ass1}, \ref{ass2} and
\ref{ass3} are satisfied. Define $j_{1}=j_{1}\left(  n\right)  $ as the
largest integer such that $2^{j_{1}}\leq n^{\frac{1}{2s+1}}$. Then, with
probability tending to one as $n\rightarrow+\infty$, the information
projection estimator (\ref{EmpBeta}) exists and satisfies:%
\[
\Delta\left(  f;f_{j_{1}(n),\widehat{\theta}_{n}}\right)  =\mathcal{O}%
_{p}\left(  n^{-\frac{2s}{2s+1}}\right)  .
\]
Moreover, the convergence is uniform over the class $F_{2,2}^{s}\left(
M\right)  $ in the sense that
\[
\lim_{K\rightarrow+\infty}\lim_{n\rightarrow+\infty}\sup_{f\in F_{2,2}%
^{s}\left(  M\right)  }{\mathbb{P}}\left(  n^{\frac{2s}{2s+1}}\Delta\left(
f;f_{j_{1}(n),\widehat{\theta}_{n}}\right)  >K\right)  =0.
\]

\end{theorem}

This theorem provides the existence with probability tending to one of a
linear estimator for the spectral density $f$ given by $f_{j_{1}\left(
n\right)  ,\widehat{\theta}_{j_{1}\left(  n\right)  }}$. This estimator is
strictly positive by construction. Therefore the corresponding estimator of
the covariance function $\widehat{\rho}^{L}$ (which is obtained as the inverse
Fourier transform of $f_{j_{1}\left(  n\right)  ,\widehat{\theta}_{n}}$) is a
positive definite function by Bochner's theorem. Hence $\widehat{\rho}^{L}$ is
a covariance function.

In the related problem of density estimation from an i.d.d. sample, Koo
\cite{K} has shown that, for the Kullback-Leibler divergence, $n^{-\frac
{2s}{2s+1}}$ is the fastest rate of convergence for the problem of estimating
a density $f$ such that $\log(f)$ belongs to the space $B_{2,2}^{s}(M)$. For
spectral densities belonging to a general Besov ball $B_{p,q}^{s}\left(
M\right)  $, Newman \cite{N} has also shown that $n^{-\frac{2s}{2s+1}}$ is an
optimal rate of convergence for the $L_{2}$ risk. For the Kullback-Leibler
divergence, we conjecture that $n^{-\frac{2s}{2s+1}}$ is the minimax rate of
convergence for spectral densities belonging to $F_{2,2}^{s}(M)$.

However, the result obtained in the above theorem is nonadaptive because the
selection of $j_{1}\left(  n\right)  $ depends on the unknown smoothness $s$
of $f$. Moreover, the result is only suited for smooth functions (as
$F_{2,2}^{s}(M)$ corresponds to a Sobolev space of order $s$) and does not
attain an optimal rate of convergence when for example $g=\log(f)$ has
singularities. We therefore propose in the next section an adaptive estimator
derived by applying an appropriate nonlinear thresholding procedure.

\subsection{Adaptive estimation}

\subsubsection{The bound on $f$ is known}

In adaptive estimation, we need to define an appropriate thresholding rule for
the wavelet coefficients of the periodogram. This threshold is level-dependent
and in this paper will take the form
\begin{equation}
\xi=\xi_{j,n}=2\left[  2\left\Vert f\right\Vert _{\infty}\left(  \sqrt
{\frac{\delta\log n}{n}}+2^{\frac{j}{2}}\left\Vert \psi\right\Vert _{\infty
}\frac{\delta\log n}{n}\right)  +\frac{C_{\ast}}{\sqrt{n}}\right]  ,
\label{threshold}%
\end{equation}
where $\delta\geq0$ is a tuning parameter whose choice will be discussed later
on and $C_{\ast}=\sqrt{\frac{C_{2}+39C_{1}^{2}}{4\pi^{2}}}$. The following
theorem states that the relative entropy between the true $f$ and its
nonlinear estimator achieves in probability the conjectured optimal rate of
convergence up to a logarithmic factor over a wide range of Besov balls.

\begin{theorem}
\label{Teo5.1enAB} Assume that $f\in F_{p,q}^{s}\left(  M\right)  $ with
$s>\frac{1}{2}+\frac{1}{p}$ and $1\leq p\leq2$. Suppose also that Assumptions
1, 2, 3 hold. For any $n>1$, define $j_{0}=j_{0}\left(  n\right)  $ to be the
integer such that $2^{j_{0}}\geq\log n\geq2^{j_{0}-1}$, and $j_{1}%
=j_{1}\left(  n\right)  $ to be the integer such that $2^{j_{1}}\geq\frac
{n}{\log n}\geq2^{j_{1}-1}$. For $\delta\geq6$, take the threshold $\xi
_{j,n}\ $as in (\ref{threshold}). Then, the thresholding estimator
(\ref{ThresBeta}) exists with probability tending to one when $n\rightarrow
+\infty$ and satisfies:%
\[
\Delta\left(  f;f_{j_{0}\left(  n\right)  ,j_{1}\left(  n\right)
,\widehat{\theta}_{n},\xi_{j,n}}^{HT}\right)  =\mathcal{O}_{p}\left(  \left(
\frac{n}{\log n}\right)  ^{-\frac{2s}{2s+1}}\right)  .
\]

\end{theorem}

Note that the choices of $j_{0}$, $j_{1}$ and $\xi_{j,n}$ are independent of
the parameter $s$; hence the estimator $f_{j_{0}\left(  n\right)
,j_{1}\left(  n\right)  ,\widehat{\theta}_{n},\xi_{j,n}}^{HT}$ is an adaptive
estimator which attains in probability what we claim is the optimal rate of
convergence, up to a logarithmic factor. In particular, $f_{j_{0}\left(
n\right)  ,j_{1}\left(  n\right)  ,\widehat{\theta}_{n},\xi_{j,n}}^{HT}$ is
adaptive on $F_{2,2}^{s}\left(  M\right)  $. This theorem provides the
existence with probability tending to one of a nonlinear estimator for the
spectral density. This estimator is strictly positive by construction.
Therefore the corresponding estimator of the covariance function
$\widehat{\rho}^{NL}$ (which is obtained as the inverse Fourier transform of
$f_{j_{0}\left(  n\right)  ,j_{1}\left(  n\right)  ,\widehat{\theta}_{n}%
,\xi_{j,n}}^{HT}$) is a positive definite function by Bochner theorem. Hence
$\widehat{\rho}^{NL}$ is a covariance function.

\subsubsection{Estimating the bound on $f$}

Although the results of Theorem \ref{Teo5.1enAB} are certainly of some
theoretical interest, they are not helpful for practical applications. The
(deterministic) threshold $\xi_{j,n}$ depends on the unknown quantities
$\left\Vert f\right\Vert _{\infty}$ and $C_{\ast}:=C\left(  C_{1}%
,C_{2}\right)  $, where $C_{1}$ and $C_{2}$ are unknown constants. To make the
method applicable, it is necessary to find some completely data-driven rule
for the threshold, which works well over a range as wide as possible of
smoothness classes. In this subsection, we give an extension that leads to
consider a random threshold which no longer depends on the bound on $f$
neither on $C_{\ast}$. For this let us consider the dyadic partitions of
$\left[  0,1\right]  $ given by $\mathcal{I}_{n}=\left\{  \left(  j/2^{J_{n}%
},(j+1)/2^{J_{n}}\right)  ,\;j=0,...,2^{J_{n}}-1\right\}  $. Given some
positive integer $r$, we define $\mathcal{P}_{n}$ as the space of piecewise
polynomials of degree $r$ on the dyadic partition $\mathcal{I}_{n}$ of step
$2^{-J_{n}}$. The dimension of $\mathcal{P}_{n}$ depends on $n$ and is denoted
by $N_{n}$. Note that $N_{n}=\left(  r+1\right)  2^{J_{n}}$. This family is
regular in the sense that the partition $\mathcal{I}_{n}$ has equispaced knots.

An estimator of $\left\Vert f\right\Vert _{\infty}$ is constructed as proposed
by Birg\'{e} and Massart \cite{B-M} in the following way. We take the infinite
norm of $\widehat{f}_{n}$, where $\widehat{f}_{n}$ denotes the (empirical)
orthogonal projection of the periodogram $I_{n}$ on $\mathcal{P}_{n}$. We
denote by $f_{n}$ the $L_{2}$-orthogonal projection of $f$ on the same space.
Then the following theorem holds.

\begin{theorem}
\label{Teo2enComte} Assume that $f\in F_{p,q}^{s}\left(  M\right)  $ with
$s>\frac{1}{2}+\frac{1}{p}$ and $1\leq p\leq2$. Suppose also that Assumptions
1, 2 and 3 hold. For any $n>1$, let $j_{0}=j_{0}\left(  n\right)  $ be the
integer such that $2^{j_{0}}\geq\log n\geq2^{j_{0}-1}$, and let $j_{1}%
=j_{1}\left(  n\right)  $ be the integer such that $2^{j_{1}}\geq\frac{n}{\log
n}\geq2^{j_{1}-1}$. Take the constants $\delta=6$ and $b\in\left[  \frac{3}%
{4},1\right)  $, and define the threshold
\begin{equation}
\widehat{\xi}_{j,n}=2\left[  2\left\Vert \widehat{f}_{n}\right\Vert _{\infty
}\left(  \sqrt{\frac{\delta}{\left(  1-b\right)  ^{2}}\frac{\log n}{n}%
}+2^{\frac{j}{2}}\left\Vert \psi\right\Vert _{\infty}\frac{\delta}{\left(
1-b\right)  ^{2}}\frac{\log n}{n}\right)  +\sqrt{\frac{\log n}{n}}\right]  .
\label{EstThreshold}%
\end{equation}
Then, if $\left\Vert f-f_{n}\right\Vert _{\infty}\leq\frac{1}{4}\left\Vert
f\right\Vert _{\infty}$ and $N_{n}\leq\frac{\kappa}{\left(  r+1\right)  ^{2}%
}\frac{n}{\log n}$, where $\kappa$ is a numerical constant and $r$ is the
degree of the polynomials, the thresholding estimator (\ref{ThresBeta}) exists
with probability tending to one as $n\rightarrow+\infty$ and satisfies%
\[
\Delta\left(  f;f_{j_{0}\left(  n\right)  ,j_{1}\left(  n\right)
,\widehat{\theta}_{n},\widehat{\xi}_{j,n}}^{HT}\right)  =\mathcal{O}%
_{p}\left(  \left(  \frac{n}{\log n}\right)  ^{-\frac{2s}{2s+1}}\right)  .
\]

\end{theorem}

Note that, we finally obtain a fully tractable estimator of $f$ which reaches
the optimal rate of convergence without prior knowledge of the regularity of
the spectral density, but also which gets rise to a real covariance estimator. \vskip.1in

\begin{remark}
We point out that, in Comte \cite{Comte} the condition $\left\Vert
f-f_{n}\right\Vert _{\infty}\leq\frac{1}{4}\left\Vert f\right\Vert _{\infty}$
is assumed. Under some regularity conditions on $f$, results from
approximation theory entails that this condition is met. Indeed for $f\in
B_{p,\infty}^{s}$, with $s>\frac{1}{p}$, we know from DeVore and Lorentz
\cite{D-L} that%
\[
\left\Vert f-f_{n}\right\Vert _{\infty}\leq C\left(  s\right)  \left\vert
f\right\vert _{s,p}N_{n}^{-\left(  s-\frac{1}{p}\right)  },
\]
with $\left\vert f\right\vert _{s,p}=\underset{y>0}{\sup}y^{-s}w_{d}\left(
f,y\right)  _{p}<+\infty$, where $w_{d}\left(  f,y\right)  _{p}$ is the
modulus of smoothness and $d=\left[  s\right]  +1$. Therefore $\left\Vert
f-f_{n}\right\Vert _{\infty}\leq\frac{1}{4}\left\Vert f\right\Vert _{\infty}$
if $N_{n}\geq\left(  4C\left(  s\right)  \frac{\left\vert f\right\vert _{s,p}%
}{\left\Vert f\right\Vert _{\infty}}\right)  ^{\frac{1}{s-\frac{1}{p}}%
}:=C\left(  f,s,p\right)  $, where $C\left(  f,s,p\right)  $ is a constant
depending on $f$, $s$ and $p$.
\end{remark}

\section{Numerical experiments}

\begin{figure}
[th]
\begin{center}
\includegraphics[
height=4.4201in,
width=7.0618in
]%
{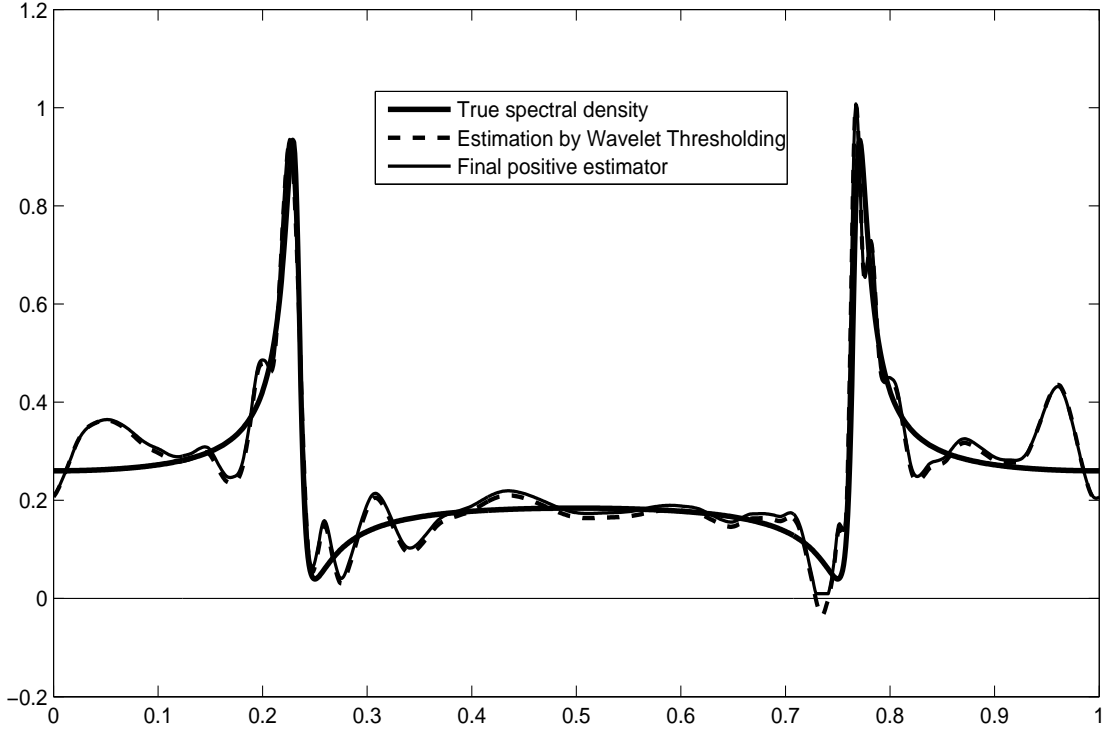}%
\caption{True spectral density $f$, wavelet thresholding estimator
$f_{j_{0},j_{1},\widehat{\xi}_{j,n}}^{HT}$ and final positive estimator
$f_{j_{0},j_{1},\widehat{\theta}_{n},\widehat{\xi}_{j,n}}^{HT}$.}%
\label{Fig1}%
\end{center}
\end{figure}

In this section we present some numerical experiments which support the claims
made in the theoretical part of this paper. The programs for our simulations
were implemented using the MATLAB programming environment. We simulate a time
series which is a superposition of an ARMA(2,2) process and a Gaussian white
noise:%
\begin{equation}
X_{t}=Y_{t}+c_{o}Z_{t}, \label{model}%
\end{equation}
where $Y_{t}+a_{1}Y_{t-1}+a_{2}Y_{t-2}=b_{0}\varepsilon_{t}+b_{1}%
\varepsilon_{t-1}+b_{2}\varepsilon_{t-2}$, and $\left\{  \varepsilon
_{t}\right\}  $, $\left\{  Z_{t}\right\}  $ are independent Gaussian white
noise processes with unit variance. The constants were chosen as $a_{1}=0.2$,
$a_{2}=0.9$, $b_{0}=1$, $b_{1}=0$, $b_{2}=1$ and $c_{0}=0.5$. We generated a
sample of size $n=1024$ according to (\ref{model}). The spectral density $f$
of $\left(  X_{t}\right)  $ is shown in Figure \ref{Fig1}. It has two
moderately sharp peaks and is smooth in the rest of the domain.

Starting from the periodogram we considered the Symmlet 8 basis, i.e. the
least asymmetric, compactly supported wavelets which are described in
Daubechies \cite{Daub}. We choose $j_{0}$ and $j_{1}$ as in the hypothesis of
Theorem \ref{Teo2enComte} and left the coefficients assigned to the father
wavelets unthresholded. Hard thresholding is performed using the threshold
$\widehat{\xi}_{j,n}$ as in (\ref{EstThreshold}) for the levels $j=j_{0}%
,...,j_{1}$, and the empirical coefficients from the higher resolution scales
$j>j_{1}$ are set to zero. This gives the estimate%
\begin{equation}
f_{j_{0},j_{1},\xi_{j,n}}^{HT}=\sum\limits_{k=0}^{2^{j_{0}}-1}\widehat
{a}_{j_{0},k}\phi_{j_{0},k}+\sum\limits_{j=j_{0}}^{j_{1}}\sum\limits_{k=0}%
^{2^{j}-1}\widehat{b}_{j,k}I\left(  \left\vert \widehat{b}_{j,k}\right\vert
>\xi_{j,n}\right)  \psi_{j,k}, \label{UWE}%
\end{equation}
which is obtained by simply thresholding the wavelet coefficients (\ref{WCE})
of the periodogram. Note that such an estimator is not guaranteed to be
strictly positive in the interval $\left[  0,1\right]  $. However, we use it
to built our strictly positive estimator $f_{j_{0},j_{1},\widehat{\theta}%
_{n},\widehat{\xi}_{j,n}}^{HT}$(see (\ref{ThresBeta}) to recall its
definition). We want to find $\widehat{\theta}_{n}$ such that%
\[
\left\langle f_{j_{0},j_{1},\widehat{\theta}_{n},\widehat{\xi}_{j,n}}%
^{HT},\psi_{j,k}\right\rangle =\delta_{\widehat{\xi}_{j,n}}\left(
\widehat{\beta}_{j,k}\right)  \mbox{ for all }\left(  j,k\right)  \in
\Lambda_{j_{1}}%
\]
For this, we take%
\[
\widehat{\theta}_{n}=\underset{\theta\in\mathbb{R}^{\#\Lambda_{j_{1}}}}%
{\arg\min}\sum\limits_{(j,k)\in\Lambda_{j_{1}}}\left(  \left\langle
f_{j_{0},j_{1},\theta},\psi_{j,k}\right\rangle -\delta_{\widehat{\xi}_{j,n}%
}\left(  \widehat{\beta}_{j,k}\right)  \right)  ^{2},
\]
where $f_{j_{0},j_{1},\theta}\left(  .\right)  =\exp\left(  \sum
\limits_{(j,k)\in\Lambda_{j_{1}}}\theta_{j,k}\psi_{j,k}\left(  .\right)
\right)  \in\mathcal{E}_{j_{1}}$ and $\mathcal{E}_{j_{1}}$ is the family
(\ref{ExpFamily}). To solve this optimization problem we used a gradient
descent method with an adaptive step, taking as initial value%
\[
\theta_{0}=\left\langle \log\left(  \left(  f_{j_{0},j_{1},\widehat{\xi}%
_{j,n}}^{HT}\right)  _{+}\right)  ,\psi_{j,k}\right\rangle ,
\]
where $\left(  f_{j_{0},j_{1},\widehat{\xi}_{j,n}}^{HT}\left(  \omega\right)
\right)  _{+}$ $:=\max\left(  f_{j_{0},j_{1},\widehat{\xi}_{j,n}}^{HT}\left(
\omega\right)  ,\eta\right)  $ for all $\omega\in\left[  0,1\right]  $ and
$\eta>0$ is a small constant.

In Figure \ref{Fig1} we display the unconstrained estimator $f_{j_{0}%
,j_{1},\xi_{j,n}}^{HT}$ as in (\ref{UWE}), obtained by thresholding of the
wavelet coefficients of the periodogram, together with the estimator
$f_{j_{0},j_{1},\widehat{\theta}_{n},\widehat{\xi}_{j,n}}^{HT}$, which is
strictly positive by construction. Note that these wavelet estimators capture
well the peaks and look fairly good on the smooth part too.%
\begin{figure}
[th]
\begin{center}
\includegraphics[
trim=0.000000in 0.000000in -0.068807in 0.136964in,
height=4.4101in,
width=7.0568in
]%
{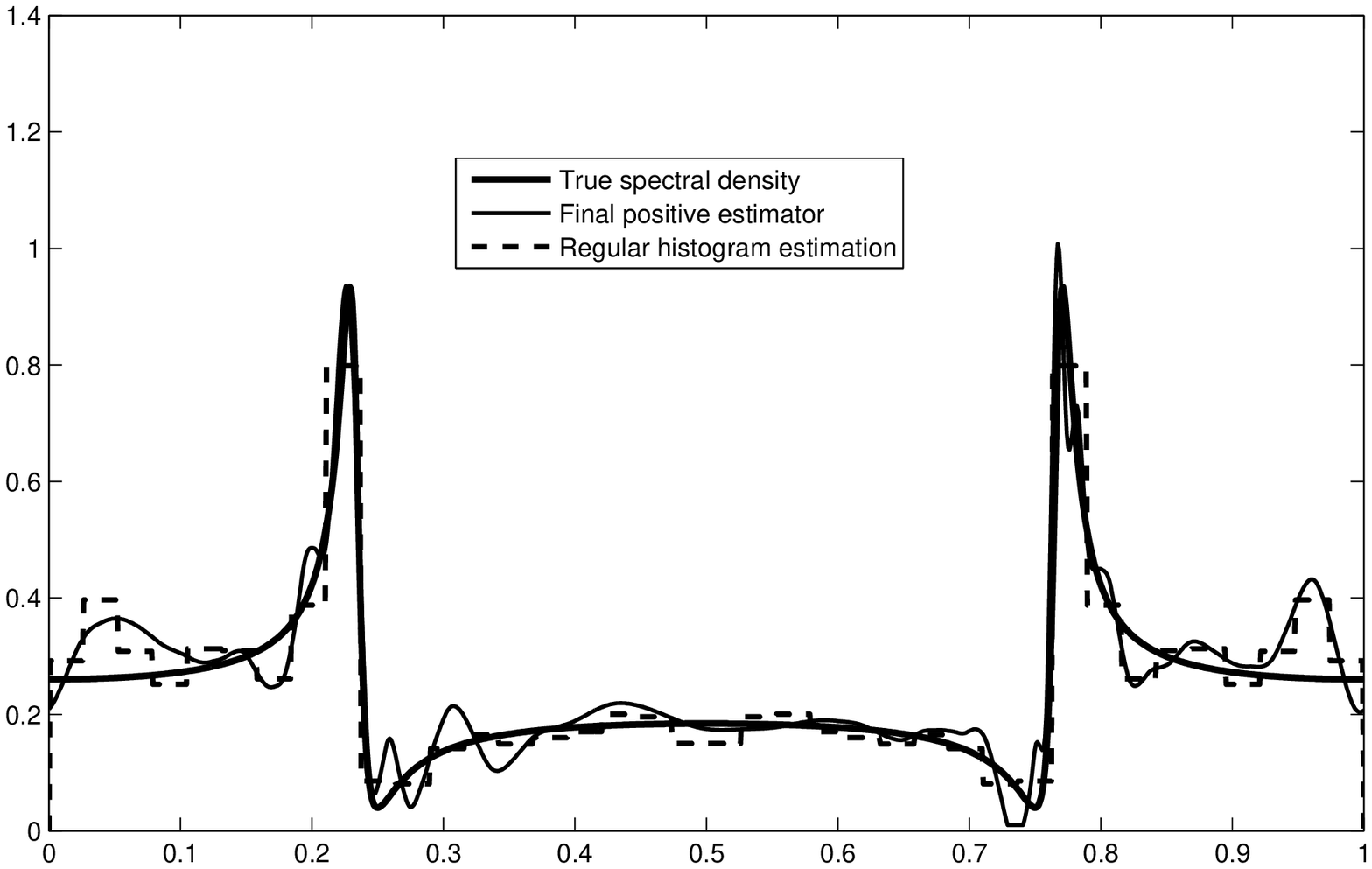}%
\caption{True spectral density $f$, final positive estimator $f_{j_{0}%
,j_{1},\widehat{\theta}_{n},\widehat{\xi}_{j,n}}^{HT}$ and estimator via model
selection using regular histograms.}%
\label{Fig2}%
\end{center}
\end{figure}

We compared our method with the spectral density estimator proposed by Comte
\cite{Comte}, which is based on a model selection procedure. As an example, in
Comte \cite{Comte}, the author study the behavior of such estimators using a
collection of nested models $\left(  S_{m}\right)  $, with $m=1,...,100$,
where $S_{m}$ is the space of piecewise constant functions, generated by a
histogram basis on $[0,1]$ of dimension $m$ with equispaced knots (see Comte
\cite{Comte} for further details). In Figure \ref{Fig2} we show the result of
this comparison. Note that our method better captures the peaks of the true
spectral density.

\section*{Acknowledgments}

This work was supported in part by Egide, under the Program of Eiffel
excellency Phd grants, as well as by the BDI CNRS grant.

\section{Appendix}

Throughout all the proofs, $C$ denotes a generic constant whose value may
change from line to line.

\subsection{Technical results for the empirical estimators of the wavelet
coefficients}

\begin{lemma}
\label{LemmaBoundBiasVarBetahat} Let $n\geq1$, $\beta_{j,k}:=\left\langle
f,\psi_{j,k}\right\rangle $ and $\widehat{\beta}_{j,k}:=\left\langle
I_{n},\psi_{j,k}\right\rangle $ for $j\geq j_{0}-1$ and $0\leq k\leq2^{j}-1$.
Suppose that Assumptions \ref{ass1}, \ref{ass2} and \ref{ass3} hold. Then,
$Bias^{2}\left(  \widehat{\beta}_{j,k}\right)  :=\left(  \mathbb{E}\left(
\widehat{\beta}_{j,k}\right)  -\beta_{j,k}\right)  ^{2}\leq\frac{C_{\ast}^{2}%
}{n}$ with $C_{\ast}=\sqrt{\frac{C_{2}+39C_{1}^{2}}{4\pi^{2}}}$, and
$Var\left(  \widehat{\beta}_{j,k}\right)  :=\mathbb{E}\left(  \widehat{\beta
}_{j,k}-\mathbb{E}\left(  \widehat{\beta}_{j,k}\right)  \right)  ^{2}\leq
\frac{C}{n}$ for some constant $C>0$. Moreover, there exists a constant
$M_{2}>0$ such that for all $f\in F_{p,q}^{s}(M)$ with $s>\frac{1}{p}$ and
$1\leq p\leq2$,
\[
{\mathbb{E}}\left(  \widehat{\beta}_{j,k}-\beta_{j,k}\right)  ^{2}%
=Bias^{2}\left(  \widehat{\beta}_{j,k}\right)  +Var\left(  \widehat{\beta
}_{j,k}\right)  \leq\frac{M_{2}}{n}.
\]

\end{lemma}

\begin{proof}
Note that $Bias^{2}\left(  \widehat{\beta}_{j,k}\right)  \leq\left\Vert
f-\mathbb{E}\left(  I_{n}\right)  \right\Vert _{L_{2}}^{2}.$ Using Proposition
1 in Comte \cite{Comte}, Assumptions \ref{ass1} and \ref{ass2} imply that
$\left\Vert f-\mathbb{E}\left(  I_{n}\right)  \right\Vert _{L_{2}}^{2}%
\leq\frac{C_{2}+39C_{1}^{2}}{4\pi^{2}n},$ which gives the result for the bias
term. To bound the variance term, remark that
\[
Var\left(  \widehat{\beta}_{j,k}\right)  =\mathbb{E}\left\langle
I_{n}-\mathbb{E}\left(  I_{n}\right)  ,\psi_{j,k}\right\rangle ^{2}%
\leq\mathbb{E}\Vert I_{n}-\mathbb{E}\left(  I_{n}\right)  \Vert_{L_{2}}%
^{2}\Vert\psi_{j,k}\Vert_{L_{2}}^{2}=\int_{0}^{1}\mathbb{E}|I_{n}%
(\omega)-\mathbb{E}\left(  I_{n}(\omega)\right)  |^{2}d\omega.
\]
Then, under Assumptions \ref{ass1} and \ref{ass2}, it follows that there
exists an absolute constant $C>0$ such that for all $\omega\in\lbrack0,1]$,
$\mathbb{E}|I_{n}(\omega)-\mathbb{E}\left(  I_{n}(\omega)\right)  |^{2}%
\leq\frac{C}{n}$. To complete the proof it remains to remark that Assumption
\ref{ass3} implies that these bounds for the bias and the variance hold
uniformly over $F_{p,q}^{s}(M)$.
\end{proof}

\begin{lemma}
\label{lem:probadeviation} Let $n\geq1$, $b_{j,k}:=\left\langle f,\psi
_{j,k}\right\rangle $ and $\widehat{b}_{j,k}:=\left\langle I_{n},\psi
_{j,k}\right\rangle $ for $j\geq j_{0}$ and $0\leq k\leq2^{j}-1$. Suppose that
Assumptions \ref{ass1} and \ref{ass2} hold. Then for any $x>0$,
\[
{\mathbb{P}}\left(  |\widehat{b}_{j,k}-b_{j,k}|>2\Vert f\Vert_{\infty}\left(
\sqrt{\frac{x}{n}}+2^{j/2}\Vert\psi\Vert_{\infty}\frac{x}{n}\right)
+\frac{C_{\ast}}{\sqrt{n}}\right)  \leq2e^{-x},
\]
where $C_{\ast}=\sqrt{\frac{C_{2}+39C_{1}^{2}}{4\pi^{2}}}$.
\end{lemma}

\begin{proof}
Note that
\[
\widehat{b}_{j,k}=\frac{1}{2\pi n}\sum\limits_{t=1}^{n}\sum\limits_{t^{\prime
}=1}^{n}\left(  X_{t}-\overline{X}\right)  \left(  X_{t^{\prime}}-\overline
{X}\right)  ^{\ast}\int\limits_{0}^{1}e^{i2\pi\omega\left(  t-t^{\prime
}\right)  }\psi_{j,k}\left(  \omega\right)  d\omega=\frac{1}{2\pi n}X^{^{T}%
}T_{n}(\psi_{j,k})X^{\ast},
\]
where $X=\left(  X_{1}-\overline{X},...,X_{n}-\overline{X}\right)  ^{T}$,
$X^{T}$ denotes the transpose of $X$ and $T_{n}(\psi_{j,k})$ is the Toeplitz
matrix with entries $\left[  T_{n}(\psi_{j,k})\right]  _{t,t^{\prime}}%
=\int\limits_{0}^{1}e^{i2\pi\omega\left(  t-t^{\prime}\right)  }\psi
_{j,k}\left(  \omega\right)  d\omega$, $1\leq t,t^{\prime}\leq n$. We can
assume without loss of generality that $E\left(  X_{t}\right)  =0$, and then
under under Assumptions \ref{ass1} and \ref{ass2}, $X$ is a centered Gaussian
vector in ${\mathbb{R}}^{n}$ with covariance matrix $\Sigma=T_{n}(f)$. Using
the decomposition $X=\Sigma^{\frac{1}{2}}\varepsilon$, where $\varepsilon\sim
N\left(  0,I_{n}\right)  $, it follows that $\widehat{b}_{j,k}=\frac{1}{2\pi
n}\varepsilon^{T}A_{j,k}\varepsilon,$ with $A_{j,k}=\Sigma^{\frac{1}{2}}%
T_{n}(\psi_{j,k})\Sigma^{\frac{1}{2}}$. Note also that $\mathbb{E}\left(
\widehat{b}_{j,k}\right)  =\frac{1}{2\pi n}tr\left(  A_{j,k}\right)  $, where
$tr\left(  A\right)  $ denotes the trace of a matrix $A$.

Now let $s_{1},\ldots,s_{n}$ be the eigenvalues of the Hermitian matrix
$A_{j,k}$ with $|s_{1}|\geq|s_{2}|\geq\ldots\geq|s_{n}|$ and let $Z=2\pi
n\left(  \widehat{b}_{j,k}-\mathbb{E}\left(  \widehat{b}_{j,k}\right)
\right)  =\varepsilon^{T}A_{j,k}\varepsilon-tr\left(  A_{j,k}\right)  $. Then,
for $0<\lambda<(2|s_{1}|)^{-1}$ one has that
\begin{align*}
\log\left(  {\mathbb{E}}\left(  e^{\lambda Z}\right)  \right)   &  =\sum
_{i=1}^{n}-\lambda s_{i}-\frac{1}{2}\log\left(  1-2\lambda s_{i}\right)
=\sum_{i=1}^{n}\sum_{\ell=2}^{+\infty}\frac{1}{2\ell}(2s_{i}\lambda)^{\ell
}\leq\sum_{i=1}^{n}\sum_{\ell=2}^{+\infty}\frac{1}{2\ell}(2|s_{i}%
|\lambda)^{\ell}\\
&  \leq\sum_{i=1}^{n}-\lambda|s_{i}|-\frac{1}{2}\log\left(  1-2\lambda
|s_{i}|\right)  ,
\end{align*}
where we have used the fact that $-\log(1-x)=\sum_{\ell=1}^{+\infty}%
\frac{x^{\ell}}{\ell}$ for $x<1$. Then using the inequality $-u-\frac{1}%
{2}\log\left(  1-2u\right)  \leq\frac{u^{2}}{1-2u}$ that holds for all
$0<u<\frac{1}{2}$, the above inequality implies that
\[
\log\left(  {\mathbb{E}}\left(  e^{\lambda Z}\right)  \right)  \leq\sum
_{i=1}^{n}\frac{\lambda^{2}|s_{i}|^{2}}{1-2\lambda|s_{i}|}\leq\frac
{\lambda^{2}\Vert s\Vert^{2}}{1-2\lambda|s_{1}|},
\]
where $\Vert s\Vert^{2}=\sum_{i=1}^{n}|s_{i}|^{2}$. Arguing as in Birg\'{e}
and Massart \cite{B-M}, the above inequality implies that for any $x>0$,
${\mathbb{P}}(|Z|>2|s_{1}|x+2\Vert s\Vert\sqrt{x})\leq2e^{-x}$, which
 implies
\begin{equation}
{\mathbb{P}}\left(  |\widehat{b}_{j,k}-\mathbb{E}\left(  \widehat{b}%
_{j,k}\right)  |>2|s_{1}|\frac{x}{n}+2\frac{\Vert s\Vert}{n}\sqrt{x}\right)
\leq2e^{-x}.\label{eq:Zx}%
\end{equation}
Let $\tau\left(  A\right)  $ denotes the spectral radius of a matrix $A$. For
the Toeplitz matrices $\Sigma=T_{n}(f)$ and $T_{n}(\psi_{j,k})$ one has that
$\tau\left(  \Sigma\right)  \leq\Vert f\Vert_{\infty}$ and $\tau\left(
T_{n}(\psi_{j,k})\right)  \leq\Vert\psi_{j,k}\Vert_{\infty}=2^{j/2}\Vert
\psi\Vert_{\infty}$. These inequalities imply that
\begin{equation}
|s_{1}|=\tau\left(  \Sigma^{\frac{1}{2}}T_{n}(\psi_{j,k})\Sigma^{\frac{1}{2}%
}\right)  \leq\tau\left(  \Sigma\right)  \tau\left(  T_{n}(\psi_{j,k})\right)
\leq\Vert f\Vert_{\infty}2^{j/2}\Vert\psi\Vert_{\infty}.\label{eq:s1}%
\end{equation}
Let $\lambda_{i}$, $i=1.,,,.n$, be the eigenvalues of $T_{n}(\psi_{j,k})$.
From Lemma 3.1 in Davies \cite{D}, we have that
\[
\underset{n\rightarrow+\infty}{\lim\sup}\frac{1}{n}tr\left(  T_{n}(\psi
_{j,k})^{2}\right)  =\underset{n\rightarrow+\infty}{\lim\sup}\frac{1}{n}%
\sum\limits_{i=1}^{n}\lambda_{i}^{2}=\int\limits_{0}^{1}\psi_{j,k}^{2}\left(
\omega\right)  d\omega=1,
\]
which implies that
\begin{equation}
\Vert s\Vert^{2}=\sum_{i=1}^{n}|s_{i}|^{2}=tr\left(  A_{j,k}^{2}\right)
=tr\left(  \left(  \Sigma T_{n}(\psi_{j,k})\right)  ^{2}\right)  \leq
\tau\left(  \Sigma\right)  ^{2}tr\left(  T_{n}(\psi_{j,k})^{2}\right)
\leq\Vert f\Vert_{\infty}^{2}n,\label{eq:s}%
\end{equation}
where we have used the inequality $tr\left(  (AB)^{2}\right)  \leq\tau\left(
A\right)  ^{2}tr\left(  B^{2}\right)  $ that holds for any pair of Hermitian
matrices $A,B$. Combining (\ref{eq:Zx}), (\ref{eq:s1}) and (\ref{eq:s}), we
finally obtain that for any $x>0$
\begin{equation}
{\mathbb{P}}\left(  |\widehat{b}_{j,k}-\mathbb{E}\left(  \widehat{b}%
_{j,k}\right)  |>2\Vert f\Vert_{\infty}\left(  \sqrt{\frac{x}{n}}+2^{j/2}%
\Vert\psi\Vert_{\infty}\frac{x}{n}\right)  \right)  \leq2e^{-x}.\label{eq:Zx2}%
\end{equation}
Now, let $\xi_{j,n}=2\Vert f\Vert_{\infty}\left(  \sqrt{\frac{x}{n}}%
+2^{j/2}\Vert\psi\Vert_{\infty}\frac{x}{n}\right)  +\frac{C_{\ast}}{\sqrt{n}}%
$, and note that%
\[
{\mathbb{P}}\left(  \left\vert \widehat{b}_{j,k}-b_{j,k}\right\vert >\xi
_{j,n}\,\right)  \leq{\mathbb{P}}\left(  \left\vert \widehat{b}_{j,k}%
-\mathbb{E}\left(  \widehat{b}_{j,k}\right)  \right\vert >\xi_{j,n}%
\,-\left\vert \mathbb{E}\left(  \widehat{b}_{j,k}\right)  -b_{j,k}\right\vert
\right)  ,
\]
By Lemma \ref{LemmaBoundBiasVarBetahat}, one has that $\left\vert
\mathbb{E}\left(  \widehat{b}_{j,k}\right)  -b_{j,k}\right\vert \leq
\frac{C_{\ast}}{\sqrt{n}}$, and thus $\xi_{j,n}\,-\left\vert \mathbb{E}\left(
\widehat{b}_{j,k}\right)  -b_{j,k}\right\vert \geq\xi_{j,n}-\frac{C_{\ast}%
}{\sqrt{n}}$ which implies using (\ref{eq:Zx2}) that
\[
{\mathbb{P}}\left(  \left\vert \widehat{b}_{j,k}-b_{j,k}\right\vert >\xi
_{j,n}\,\right)  \leq{\mathbb{P}}\left(  \left\vert \widehat{b}_{j,k}%
-\mathbb{E}\left(  \widehat{b}_{j,k}\right)  \right\vert >\xi_{j,n}%
-\frac{C_{\ast}}{\sqrt{n}}\right)  \leq2e^{-x},
\]
which completes the proof of Lemma \ref{lem:probadeviation}.
\end{proof}

\begin{lemma}
\label{LemmaBoundBiasVarThresBetahat} Assume that $f\in F_{p,q}^{s}\left(
M\right)  $ with $s>\frac{1}{2}+\frac{1}{p}$ and $1\leq p\leq2$. Suppose that
Assumptions \ref{ass1}, \ref{ass2} and \ref{ass3} hold. For any $n>1$, define
$j_{0}=j_{0}\left(  n\right)  $ to be the integer such that $2^{j_{0}}>\log
n\geq2^{j_{0}-1}$, and $j_{1}=j_{1}\left(  n\right)  $ to be the integer such
that $2^{j_{1}}\geq\frac{n}{\log n}\geq2^{j_{1}-1}$. For $\delta\geq6$, take
the threshold $\xi_{j,n}=2\left[  2\left\Vert f\right\Vert _{\infty}\left(
\sqrt{\frac{\delta\log n}{n}}+2^{\frac{j}{2}}\left\Vert \psi\right\Vert
_{\infty}\frac{\delta\log n}{n}\right)  +\frac{C_{\ast}}{\sqrt{n}}\right]  $
as in (\ref{threshold}), where $C_{\ast}=\sqrt{\frac{C_{2}+39C_{1}^{2}}%
{4\pi^{2}}}$. Let $\beta_{j,k}:=\left\langle f,\psi_{j,k}\right\rangle $ and
$\widehat{\beta}_{\xi_{j,n},\left(  j,k\right)  }:=\delta_{\xi_{j,n}}\left(
\widehat{\beta}_{j,k}\right)  $ with $\left(  j,k\right)  \in$ $\Lambda
_{j_{1}}$ as in (\ref{ThresBeta}). Take $\beta=\left(  \beta_{j,k}\right)
_{\left(  j,k\right)  \in\Lambda_{j_{1}}}$ and $\widehat{\beta}_{\xi_{j,n}%
}=\left(  \widehat{\beta}_{\xi_{j,n},\left(  j,k\right)  }\right)  _{\left(
j,k\right)  \in\Lambda_{j_{1}}}$. Then there exists a constant $M_{3}>0$ such
that for all sufficiently large $n$:
\[
\mathbb{E}\left\Vert \beta-\widehat{\beta}_{\xi_{j,n}}\right\Vert _{2}%
^{2}:=\mathbb{E}\left(  \sum\limits_{\left(  j,k\right)  \in\Lambda_{j_{1}}%
}\left\vert \beta_{j,k}-\delta_{\xi_{j,n}}\left(  \widehat{\beta}%
_{j,k}\right)  \right\vert ^{2}\right)  \leq M_{3}\left(  \frac{n}{\log
n}\right)  ^{-\frac{2s}{2s+1}}%
\]
uniformly over $F_{p,q}^{s}(M)$.
\end{lemma}

\begin{proof}
Taking into account that%
\begin{align}
\mathbb{E}\left\Vert \beta-\widehat{\beta}_{\xi_{j,n}}\right\Vert _{2}^{2} &
=\sum\limits_{k=0}^{2^{j_{0}}-1}\mathbb{E}\left(  a_{j_{0},k}-\widehat
{a}_{j_{0},k}\right)  ^{2}+\sum\limits_{j=j_{0}}^{j_{1}}\sum\limits_{k=0}%
^{2^{j}-1}\mathbb{E}\left[  \left(  b_{j,k}-\widehat{b}_{j,k}\right)
^{2}I\left(  \left\vert \widehat{b}_{j,k}\right\vert >\xi_{j,n}\right)
\right]  \nonumber\\
&  +\sum\limits_{j=j_{0}}^{j_{1}}\sum\limits_{k=0}^{2^{j}-1}b_{j,k}%
^{2}P\left(  \left\vert \widehat{b}_{j,k}\right\vert \leq\xi_{j,n}\right)
\nonumber\\
&  :=T_{1}+T_{2}+T_{3},\label{Sum6terms}%
\end{align}
we are interested in bounding these three terms. The bound for $T_{1}$ follows
from Lemma \ref{LemmaBoundBiasVarBetahat} and the fact that $j_{0}=\log
_{2}\left(  \log n\right)  \leq\frac{1}{2s+1}\log_{2}\left(  n\right)  $:%
\begin{equation}
T_{1}=\sum\limits_{k=0}^{2^{j_{0}}-1}\mathbb{E}\left(  a_{j_{0},k}-\widehat
{a}_{j_{0},k}\right)  ^{2}=O\left(  \frac{2^{j_{0}}}{n}\right)  \leq O\left(
n^{-\frac{2s}{2s+1}}\right)  .\label{BforT1}%
\end{equation}
To bound $T_{2}$ and $T_{3}$ we proceed as follows. Write%
\[
T_{2}=\sum\limits_{j=j_{0}}^{j_{1}}\sum\limits_{k=0}^{2^{j}-1}\mathbb{E}%
\left[  \left(  b_{j,k}-\widehat{b}_{j,k}\right)  ^{2}\left\{  I\left(
\left\vert \widehat{b}_{j,k}\right\vert >\xi_{j,n}\,,\left\vert b_{j,k}%
\right\vert >\frac{\xi_{j,n}}{2}\,\right)  +I\left(  \left\vert \widehat
{b}_{j,k}\right\vert >\xi_{j,n}\,,\left\vert b_{j,k}\right\vert \leq\frac
{\xi_{j,n}}{2}\,\right)  \right\}  \right]
\]
and
\[
T_{3}=\sum\limits_{j=j_{0}}^{j_{1}}\sum\limits_{k=0}^{2^{j}-1}b_{j,k}%
^{2}\left[  P\left(  \left\vert \widehat{b}_{j,k}\right\vert \leq\xi
_{j,n}\,,\left\vert b_{j,k}\right\vert \leq2\xi_{j,n}\right)  +P\left(
\left\vert \widehat{b}_{j,k}\right\vert \leq\xi_{j,n},\left\vert
b_{j,k}\right\vert >2\xi_{j,n}\right)  \right]  .
\]
From Hardle, Kerkyacharian, Picard and Tsybakov \cite{H-K-P-T} we get that%
\begin{align*}
T_{2}+T_{3} &  \leq\sum\limits_{j=j_{0}}^{j_{1}}\sum\limits_{k=0}^{2^{j}%
-1}\mathbb{E}\left\{  \left(  b_{j,k}-\widehat{b}_{j,k}\right)  ^{2}\right\}
I\left(  \left\vert b_{j,k}\right\vert >\frac{\xi_{j,n}}{2}\right)
+\sum\limits_{j=j_{0}}^{j_{1}}\sum\limits_{k=0}^{2^{j}-1}b_{j,k}^{2}I\left(
\left\vert b_{j,k}\right\vert \leq2\xi_{j,n}\right)  \\
&  +5\sum\limits_{j=j_{0}}^{j_{1}}\sum\limits_{k=0}^{2^{j}-1}\mathbb{E}%
\left\{  \left(  b_{j,k}-\widehat{b}_{j,k}\right)  ^{2}I\left(  \left\vert
\widehat{b}_{j,k}-b_{j,k}\right\vert >\frac{\xi_{j,n}}{2}\,\right)  \right\}
\\
&  :=T^{\prime}+T^{\prime\prime}+T^{\prime\prime\prime}.
\end{align*}
Now we bound $T^{\prime\prime\prime}$. Using Cauchy-Schwarz inequality, we
obtain%
\[
T^{\prime\prime\prime}\leq5\sum\limits_{j=j_{0}}^{j_{1}}\sum\limits_{k=0}%
^{2^{j}-1}\mathbb{E}^{\frac{1}{2}}\left[  \left(  b_{j,k}-\widehat{b}%
_{j,k}\right)  ^{4}\right]  P^{\frac{1}{2}}\left(  \left\vert \widehat
{b}_{j,k}-b_{j,k}\right\vert >\frac{\xi_{j,n}}{2}\,\right)  .
\]
By the same inequality we get $\mathbb{E}\left[  \left(  \widehat{b}%
_{j,k}-b_{j,k}\right)  ^{4}\right]  \leq\mathbb{E}\left[  \left\Vert
I_{n}-f\right\Vert _{L_{2}}^{4}\left\Vert \psi_{j,k}\right\Vert _{L_{2}}%
^{4}\right]  =O\left(  \mathbb{E}\left\Vert I_{n}-f\right\Vert _{L_{2}}%
^{4}\right)  $. It can be checked that $\mathbb{E}\left\Vert I_{n}%
-f\right\Vert _{L_{2}}^{4}\leq8\mathbb{E}\left(  \left\Vert I_{n}%
-\mathbb{E}\left(  I_{n}\right)  \right\Vert _{L_{2}}^{4}+\left\Vert
\mathbb{E}\left(  I_{n}\right)  -f\right\Vert _{L_{2}}^{4}\right)  $.
According to Comte \cite{Comte}, $\mathbb{E}\left\Vert I_{n}-\mathbb{E}\left(
I_{n}\right)  \right\Vert _{L_{2}}^{4}=O\left(  n^{2}\right)  $. From the
proof of Lemma \ref{LemmaBoundBiasVarBetahat} we get that $\left\Vert
\mathbb{E}\left(  I_{n}\right)  -f\right\Vert _{L_{2}}^{4}=O\left(  \frac
{1}{n^{2}}\right)  $. Therefore $\mathbb{E}\left\Vert I_{n}-f\right\Vert
_{L_{2}}^{4}\leq O\left(  n^{2}+\frac{1}{n^{2}}\right)  =O\left(
n^{2}\right)  $. Hence $\mathbb{E}\left[  \left(  \widehat{b}_{j,k}%
-b_{j,k}\right)  ^{4}\right]  =O\left(  \mathbb{E}\left\Vert I_{n}%
-f\right\Vert _{L_{2}}^{4}\right)  =O\left(  n^{2}\right)  $. For the bound of
$P\left(  \left\vert \widehat{b}_{j,k}-b_{j,k}\right\vert >\frac{\xi_{j,n}}%
{2}\,\right)  $ we use the result of Lemma \ref{lem:probadeviation} with
$x=\delta\log n$, where $\delta>0$ is a constant to be specified later. We
obtain%
\begin{align*}
P\left(  \left\vert \widehat{b}_{j,k}-b_{j,k}\right\vert >\frac{\xi_{j,n}}%
{2}\,\right)   &  =P\left(  \left\vert \widehat{b}_{j,k}-b_{j,k}\right\vert
>2\left\Vert f\right\Vert _{\infty}\left(  \sqrt{\frac{\delta\log n}{n}%
}+2^{\frac{j}{2}}\left\Vert \psi\right\Vert _{\infty}\frac{\delta\log n}%
{n}\right)  +\frac{C_{\ast}}{\sqrt{n}}\right)  \\
&  \leq2e^{-\delta\log n}=2n^{-\delta}.
\end{align*}
Therefore, for $\delta\geq6$, we get%
\[
T^{\prime\prime\prime}\leq5\sum\limits_{j=j_{0}}^{j_{1}}\sum\limits_{k=0}%
^{2^{j}-1}\mathbb{E}^{\frac{1}{2}}\left[  \left(  b_{j,k}-\widehat{b}%
_{j,k}\right)  ^{4}\right]  P^{\frac{1}{2}}\left(  \left\vert \widehat
{b}_{j,k}-b_{j,k}\right\vert >\frac{\xi_{j,n}}{2}\,\right)  \leq O\left(
\frac{n^{-1}}{\log n}\right)  \leq O\left(  n^{-\frac{2s}{2s+1}}\right)  .
\]

Now we follow results found in Pensky and Sapatinas \cite{P-S} to bound
$T^{\prime}$ and $T^{\prime\prime}$. Let $j_{A}$ be the integer such that
$2^{j_{A}}>\left(  \frac{n}{\log n}\right)  ^{\frac{1}{2s+1}}>2^{j_{A}-1}$
(note that given our assumptions $j_{0}\leq j_{A}\leq j_{1}$ for all
sufficiently large $n$), then $T^{\prime}$ can be partitioned as $T^{\prime
}=T_{1}^{\prime}+T_{2}^{\prime}$, where the first component is calculated over
the set of indices $j_{0}\leq j\leq j_{A}$ and the second component over
$j_{A}+1\leq j\leq j_{1}$. Hence, using Lemma \ref{LemmaBoundBiasVarBetahat}
we obtain%
\begin{equation}
T_{1}^{\prime}\leq C\sum\limits_{j=j_{0}}^{j_{A}}\frac{2^{j}}{n}=O\left(
2^{j_{A}}n^{-1}\right)  =O\left(  \left(  \frac{n}{\log n}\right)  ^{\frac
{1}{2s+1}}n^{-1}\right)  \leq O\left(  n^{-\frac{2s}{2s+1}}\right)
.\label{BforT1I}%
\end{equation}
To obtain a bound for $T_{2}^{\prime}$, we will use that if $f\in F_{p,q}%
^{s}\left(  A\right)  $, then for some constant $C$, dependent on $s$, $p$,
$q$ and $A>0$ only, we have that%
\begin{equation}
\sum\limits_{k=0}^{2^{j}-1}b_{j,k}^{2}\leq C2^{-2js^{\ast}},\label{cnd}%
\end{equation}
for $1\leq p\leq2$, where $s^{\ast}=s+\frac{1}{2}-\frac{1}{p}$. Taking into
account that $I\left(  \left\vert b_{j,k}\right\vert >\frac{\xi_{j,n}}%
{2}\right)  \leq\frac{4}{\xi_{j,n}^{2}}\left\vert b_{j,k}\right\vert ^{2}$, we
get
\begin{align*}
T_{2}^{\prime} &  \leq\frac{C}{n}\sum\limits_{j=j_{A}}^{j_{1}}\sum
\limits_{k=0}^{2^{j}-1}\frac{4}{\xi_{j,n}^{2}}\left\vert b_{j,k}\right\vert
^{2}\leq\frac{C\left(  \left\Vert f\right\Vert _{\infty}\right)  2^{-2s^{\ast
}j_{A}}}{\left(  \sqrt{\delta\log n}+\left\Vert \psi\right\Vert _{\infty
}\delta n^{\frac{-s}{2s+1}}\left(  \log n\right)  ^{\frac{4s+1}{4s+2}}\right)
^{2}}\sum\limits_{j=j_{A}}^{j_{1}}\sum\limits_{k=0}^{2^{j}-1}2^{2js^{\ast}%
}\left\vert b_{j,k}\right\vert ^{2}\\
&  \leq O\left(  2^{-2s^{\ast}j_{A}}\right)  =O\left(  \left(  \frac{n}{\log
n}\right)  ^{-\frac{2s^{\ast}}{2s+1}}\right)  ,
\end{align*}
where we used the fact that $\sqrt{\delta\log n}+\left\Vert \psi\right\Vert
_{\infty}\delta n^{\frac{-s}{2s+1}}\left(  \log n\right)  ^{\frac{4s+1}{4s+2}%
}\rightarrow+\infty$ when $n\rightarrow+\infty$. Now remark that if $p=2$ then
$s^{\ast}=s$ and thus
\begin{equation}
T_{2}^{\prime}=O\left(  \left(  \frac{n}{\log n}\right)  ^{-\frac{2s^{\ast}%
}{2s+1}}\right)  =O\left(  \left(  \frac{n}{\log n}\right)  ^{-\frac{2s}%
{2s+1}}\right)  .\label{B1forT2I}%
\end{equation}
For the case $1\leq p<2$, the repeated use of the fact that if $B,D>0$ then
$I\left(  \left\vert b_{j,k}\right\vert >B+D\right)  \leq I\left(  \left\vert
b_{j,k}\right\vert >B\right)  $, enables us to obtain that
\begin{align*}
T_{2}^{\prime} &  \leq\frac{C}{n}\sum\limits_{j=j_{A}}^{j_{1}}\sum
\limits_{k=0}^{2^{j}-1}I\left(  \left\vert b_{j,k}\right\vert >\frac{\xi
_{j,n}}{2}\right)  \leq\frac{C}{n}\sum\limits_{j=j_{A}}^{j_{1}}\sum
\limits_{k=0}^{2^{j}-1}\left\vert b_{j,k}\right\vert ^{-p}\left\vert
b_{j,k}\right\vert ^{p}I\left(  \left\vert b_{j,k}\right\vert ^{-p}<\left(
2\left\Vert f\right\Vert _{\infty}\sqrt{\delta}\sqrt{\frac{\log n}{n}}\right)
^{-p}\right)  \\
&  \leq C\sum\limits_{j=j_{A}}^{j_{1}}\sum\limits_{k=0}^{2^{j}-1}\frac{1}%
{n}\left(  2\left\Vert f\right\Vert _{\infty}\sqrt{\delta}\sqrt{\frac{\log
n}{n}}\right)  ^{-p}\left\vert b_{j,k}\right\vert ^{p}.
\end{align*}
Since $f\in F_{p,q}^{s}\left(  A\right)  $ it follows that there exists a
constant $C$ depending only on $p$, $q$, $s$ and $A$ such that%
\begin{equation}
\sum\limits_{k=0}^{2^{j}-1}\left\vert b_{j,k}\right\vert ^{p}\leq
C2^{-pjs^{\ast}},\label{cnd1}%
\end{equation}
where $s^{\ast}=s+\frac{1}{2}-\frac{1}{p}$ as before. By (\ref{cnd1}) we get%
\begin{align}
T_{2}^{\prime} &  \leq\left(  \log n\right)  C\left(  \left\Vert f\right\Vert
_{\infty},\delta,p\right)  \sum\limits_{j=j_{A}}^{j_{1}}\sum\limits_{k=0}%
^{2^{j}-1}\frac{\left(  \log n\right)  ^{-\frac{p}{2}}}{n^{1-\frac{p}{2}}%
}\left\vert b_{j,k}\right\vert ^{p}\leq C\left(  \left\Vert f\right\Vert
_{\infty},\delta,p\right)  \frac{\left(  \log n\right)  ^{1-\frac{p}{2}}%
}{n^{1-\frac{p}{2}}}\sum\limits_{j=j_{A}}^{j_{1}}C2^{-pjs^{\ast}}\nonumber\\
&  =O\left(  \frac{\left(  \log n\right)  ^{1-\frac{p}{2}}}{n^{1-\frac{p}{2}}%
}2^{-pj_{A}s^{\ast}}\right)  =O\left(  \left(  \frac{n}{\log n}\right)
^{-\frac{2s}{2s+1}}\right)  .\label{B2forT2Iaux}%
\end{align}
Hence, by (\ref{BforT1I}), (\ref{B1forT2I}) and (\ref{B2forT2Iaux}),
$T^{\prime}=O\left(  \left(  \frac{n}{\log n}\right)  ^{-\frac{2s}{2s+1}%
}\right)  $.

Now, set $j_{A}$ as before, then $T^{\prime\prime}$ can be split into
$T^{\prime\prime}=T_{1}^{\prime\prime}+T_{2}^{\prime\prime}$, where the first
component is calculated over the set of indices $j_{0}\leq j\leq j_{A}$ and
the second component over $j_{A}+1\leq j\leq j_{1}$. Then
\[
T_{1}^{\prime\prime}\leq\sum\limits_{j=j_{0}}^{j_{A}}\sum\limits_{k=0}%
^{2^{j}-1}b_{j,k}^{2}I\left(  \left\vert b_{j,k}\right\vert ^{2}\leq32\left[
4\left\Vert f\right\Vert _{\infty}^{2}\left(  \sqrt{\frac{\delta\log n}{n}%
}+2^{\frac{j}{2}}\left\Vert \psi\right\Vert _{\infty}\frac{\delta\log n}%
{n}\right)  ^{2}+\frac{C_{\ast}^{2}}{n}\right]  \right)  .
\]
Using repeatedly that $\left(  B+D\right)  ^{2}\leq2\left(  B^{2}%
+D^{2}\right)  $ for $B,D\in\mathbb{R}$, we obtain the desired bound for
$T_{1}^{\prime\prime}$:%
\begin{align}
T_{1}^{\prime\prime} &  \leq C\left(  \left\Vert f\right\Vert _{\infty
}\right)  \sum\limits_{j=j_{0}}^{j_{A}}\sum\limits_{k=0}^{2^{j}-1}\left(
\frac{\delta\log n}{n}+2^{j}\left\Vert \psi\right\Vert _{\infty}^{2}%
\frac{\delta^{2}\left(  \log n\right)  ^{2}}{n^{2}}\right)  +C\left(  C_{\ast
}\right)  \sum\limits_{j=j_{0}}^{j_{A}}\sum\limits_{k=0}^{2^{j}-1}\frac{1}%
{n}\nonumber\\
&  \leq C\left(  \left\Vert f\right\Vert _{\infty},\delta,C_{\ast}\right)
\frac{\log n}{n}2^{j_{A}}+C\left(  \left\Vert f\right\Vert _{\infty}%
,\delta,\left\Vert \psi\right\Vert _{\infty}\right)  \frac{\left(  \log
n\right)  ^{2}}{n^{2}}2^{2j_{A}}\nonumber\\
&  =O\left(  \left(  \log n\right)  ^{\frac{2s}{2s+1}}n^{-\frac{2s}{2s+1}%
}+\left(  \log n\right)  ^{\frac{4s}{2s+1}}n^{-\frac{4s}{2s+1}}\right)  \leq
O\left(  \left(  \frac{n}{\log n}\right)  ^{-\frac{2s}{2s+1}}\right)
.\label{BforT1II}%
\end{align}
To bound $T_{2}^{\prime\prime}$, note that $T_{2}^{\prime\prime}\leq
\sum\limits_{j=j_{A}+1}^{j_{1}}\sum\limits_{k=0}^{2^{j}-1}b_{j,k}^{2}=O\left(
2^{-2j_{A}s^{\ast}}\right)  =O\left(  \left(  \frac{n}{\log n}\right)
^{-\frac{2s^{\ast}}{2s+1}}\right)  $, where we have used the condition
(\ref{cnd}). Now remark that if $p=2$ then $s^{\ast}=s$ and thus
\begin{equation}
T_{2}^{\prime\prime}=O\left(  \left(  \frac{n}{\log n}\right)  ^{-\frac
{2s^{\ast}}{2s+1}}\right)  =O\left(  \left(  \frac{n}{\log n}\right)
^{-\frac{2s}{2s+1}}\right)  .\label{B1forT2II}%
\end{equation}
If $1\leq p<2$,%
\begin{align}
T_{2}^{\prime\prime} &  =\sum\limits_{j=j_{A}+1}^{j_{1}}\sum\limits_{k=0}%
^{2^{j}-1}\left\vert b_{j,k}\right\vert ^{2-p}\left\vert b_{j,k}\right\vert
^{p}I\left(  \left\vert b_{j,k}\right\vert \leq2\xi_{j,n}\,\right)
\nonumber\\
&  \leq\sum\limits_{j=j_{A}+1}^{j_{1}}\sum\limits_{k=0}^{2^{j}-1}\left(
8\left\Vert f\right\Vert _{\infty}\sqrt{\frac{\delta\log n}{n}}+8\left\Vert
f\right\Vert _{\infty}2^{\frac{j_{1}}{2}}\left\Vert \psi\right\Vert _{\infty
}\frac{\delta\log n}{n}+4\sqrt{\frac{\log n}{n}}\right)  ^{2-p}\left\vert
b_{j,k}\right\vert ^{p}\nonumber\\
&  \leq\left(  C\left(  \left\Vert f\right\Vert _{\infty},\left\Vert
\psi\right\Vert _{\infty},\delta\right)  \right)  ^{2-p}\left(  \sqrt
{\frac{\log n}{n}}\right)  ^{2-p}\sum\limits_{j=j_{A}+1}^{j_{1}}%
\sum\limits_{k=0}^{2^{j}-1}\left\vert b_{j,k}\right\vert ^{p}=O\left(  \left(
\frac{\log n}{n}\right)  ^{\frac{2-p}{2}}2^{-pj_{A}s^{\ast}}\right)
\nonumber\\
&  =O\left(  \left(  \frac{n}{\log n}\right)  ^{\frac{p}{2}-1-\frac{p\left(
s+\frac{1}{2}-\frac{1}{p}\right)  }{2s+1}}\right)  =O\left(  \left(  \frac
{n}{\log n}\right)  ^{-\frac{2s}{2s+1}}\right)  ,\label{B2forT2II}%
\end{align}
where we have used condition (\ref{cnd1}) and the fact that $C_{\ast}\leq
\sqrt{\log n}$ for $n$ sufficiently large, taking into account that the
constant $C_{\ast}:=C\left(  C_{1},C_{2}\right)  $ does not depend on $n$.
Hence, by (\ref{BforT1II}), (\ref{B1forT2II}) and (\ref{B2forT2II}),
$T^{\prime\prime}=O\left(  \left(  \frac{n}{\log n}\right)  ^{-\frac{2s}%
{2s+1}}\right)  $. Combining all terms in (\ref{Sum6terms}), we conclude that:%
\[
\mathbb{E}\left\Vert \beta-\widehat{\beta}_{\xi_{j,n}}\right\Vert _{2}%
^{2}=O\left(  \left(  \frac{n}{\log n}\right)  ^{-\frac{2s}{2s+1}}\right)  .
\]
This completes the proof.
\end{proof}

\begin{lemma}
\label{LemmaBoundBiasVarThresBetahatEpsHat} Assume that $f\in F_{p,q}%
^{s}\left(  M\right)  $ with $s>\frac{1}{2}+\frac{1}{p}$ and $1\leq p\leq2$.
Suppose that Assumptions \ref{ass1}, \ref{ass2} and \ref{ass3} hold. For any
$n>1$, define $j_{0}=j_{0}\left(  n\right)  $ to be the integer such that
$2^{j_{0}}>\log n\geq2^{j_{0}-1}$, and $j_{1}=j_{1}\left(  n\right)  $ to be
the integer such that $2^{j_{1}}\geq\frac{n}{\log n}\geq2^{j_{1}-1}$. Define
the threshold $\widehat{\xi}_{j,n}$ as in (\ref{EstThreshold}) for some
constants $\delta=6$ and $b\in\left[  \frac{3}{4},1\right)  $. Let
$\beta_{j,k}:=\left\langle f,\psi_{j,k}\right\rangle $ and $\widehat{\beta
}_{\widehat{\xi}_{j,n},\left(  j,k\right)  }:=\delta_{\widehat{\xi}_{j,n}%
}\left(  \widehat{\beta}_{j,k}\right)  $ with $\left(  j,k\right)  \in$
$\Lambda_{j_{1}}$ as in (\ref{ThresBeta}). Take $\beta=\left(  \beta
_{j,k}\right)  _{\left(  j,k\right)  \in\Lambda_{j_{1}}}$ and $\widehat{\beta
}_{\widehat{\xi}_{j,n}}=\left(  \widehat{\beta}_{\widehat{\xi}_{j,n},\left(
j,k\right)  }\right)  _{\left(  j,k\right)  \in\Lambda_{j_{1}}}$. Then, if
$\left\Vert f-f_{n}\right\Vert _{\infty}\leq\frac{1}{4}\left\Vert f\right\Vert
_{\infty}$ and $N_{n}\leq\frac{\kappa}{\left(  r+1\right)  ^{2}}\frac{n}{\log
n}$, where $\kappa$ is a numerical constant and $r$ is the degree of the
polynomials, there exists a constant $M_{4}>0$ such that for all sufficiently
large $n$:
\[
\mathbb{E}\left\Vert \beta-\widehat{\beta}_{\widehat{\xi}_{j,n}}\right\Vert
_{2}^{2}:=\mathbb{E}\left(  \sum\limits_{\left(  j,k\right)  \in\Lambda
_{j_{1}}}\left\vert \delta_{\widehat{\xi}_{j,n}}\left(  \widehat{\beta}%
_{j,k}\right)  -\beta_{j,k}\right\vert ^{2}\right)  \leq M_{4}\left(  \frac
{n}{\log n}\right)  ^{-\frac{2s}{2s+1}}%
\]
uniformly over $F_{p,q}^{s}(M)$.
\end{lemma}

\begin{proof}
Recall that $f_{n}$ is the $L_{2}$ orthogonal projection of $f$ on the space
$\mathcal{P}_{n}$ of piecewise polynomials of degree $r$ on a dyadic partition
with step $2^{-J_{n}}$. The dimension of $\mathcal{P}_{n}$ is $N_{n}=\left(
r+1\right)  2^{J_{n}}$. Let $\widehat{f}_{n}$ similarly be the orthogonal
projection of $I_{n}$ on $\mathcal{P}_{n}$. By doing analogous work as the one
done to obtain (\ref{Sum6terms}), we get that%
\begin{equation}
\mathbb{E}\left\Vert \beta-\widehat{\beta}_{\widehat{\xi}_{j,n}}\right\Vert
_{2}^{2}:=T_{1}+T_{2}+T_{3},\label{Sum4terminos}%
\end{equation}
where $T_{1}=\sum\limits_{k=0}^{2^{j_{0}}-1}\mathbb{E}\left(  a_{j_{0}%
,k}-\widehat{a}_{j_{0},k}\right)  ^{2}$ do not depend on $\widehat{\xi}_{j,n}%
$. Therefore, by (\ref{BforT1}), $T_{1}=O\left(  n^{-\frac{2s}{2s+1}}\right)
$. For $T_{2}$ and $T_{3}$ we have that%
\[
T_{2}=\sum\limits_{j=j_{0}}^{j_{1}}\sum\limits_{k=0}^{2^{j}-1}\mathbb{E}%
\left[  \left(  b_{j,k}-\widehat{b}_{j,k}\right)  ^{2}\left\{  I\left(
\left\vert \widehat{b}_{j,k}\right\vert >\widehat{\xi}_{j,n}\,,\left\vert
b_{j,k}\right\vert >\frac{\widehat{\xi}_{j,n}}{2}\,\right)  +I\left(
\left\vert \widehat{b}_{j,k}\right\vert >\widehat{\xi}_{j,n}\,,\left\vert
b_{j,k}\right\vert \leq\frac{\widehat{\xi}_{j,n}}{2}\,\right)  \right\}
\right]
\]
and%
\[
T_{3}=\sum\limits_{j=j_{0}}^{j_{1}}\sum\limits_{k=0}^{2^{j}-1}b_{j,k}%
^{2}\left[  P\left(  \left\vert \widehat{b}_{j,k}\right\vert \leq\widehat{\xi
}_{j,n}\,,\left\vert b_{j,k}\right\vert \leq2\widehat{\xi}_{j,n}\right)
+P\left(  \left\vert \widehat{b}_{j,k}\right\vert \leq\widehat{\xi}%
_{j,n},\left\vert b_{j,k}\right\vert >2\widehat{\xi}_{j,n}\right)  \right]  .
\]
Using the same decomposition as in the proof of Lemma
(\ref{LemmaBoundBiasVarThresBetahat}) we get that%
\begin{align*}
T_{2}+T_{3} &  \leq\sum\limits_{j=j_{0}}^{j_{1}}\sum\limits_{k=0}^{2^{j}%
-1}\mathbb{E}\left\{  \left(  b_{j,k}-\widehat{b}_{j,k}\right)  ^{2}I\left(
\left\vert b_{j,k}\right\vert >\frac{\widehat{\xi}_{j,n}}{2}\right)  \right\}
+\sum\limits_{j=j_{0}}^{j_{1}}\sum\limits_{k=0}^{2^{j}-1}\mathbb{E}\left\{
b_{j,k}^{2}I\left(  \left\vert b_{j,k}\right\vert \leq2\widehat{\xi}%
_{j,n}\right)  \right\}  \\
&  +5\sum\limits_{j=j_{0}}^{j_{1}}\sum\limits_{k=0}^{2^{j}-1}\mathbb{E}%
\left\{  \left(  b_{j,k}-\widehat{b}_{j,k}\right)  ^{2}I\left(  \left\vert
\widehat{b}_{j,k}-b_{j,k}\right\vert >\frac{\widehat{\xi}_{j,n}}{2}\,\right)
\right\}  :=T^{\prime}+T^{\prime\prime}+T^{\prime\prime\prime}.
\end{align*}
Now we bound $T^{\prime\prime\prime}$. Using Cauchy-Schwarz inequality, one
obtains%
\[
T^{\prime\prime\prime}\leq5\sum\limits_{j=j_{0}}^{j_{1}}\sum\limits_{k=0}%
^{2^{j}-1}\mathbb{E}^{\frac{1}{2}}\left[  \left(  b_{j,k}-\widehat{b}%
_{j,k}\right)  ^{4}\right]  P^{\frac{1}{2}}\left(  \left\vert \widehat
{b}_{j,k}-b_{j,k}\right\vert >\frac{\widehat{\xi}_{j,n}}{2}\,\right)  ,
\]
From Lemma \ref{lem:probadeviation} we have that for any $y>0$ the following
exponential inequality holds:%
\begin{equation}
P\left(  \left\vert \widehat{b}_{j,k}-b_{j,k}\right\vert >2\left\Vert
f\right\Vert _{\infty}\left(  \sqrt{\frac{y}{n}}+2^{\frac{j}{2}}\left\Vert
\psi\right\Vert _{\infty}\frac{y}{n}\right)  +\frac{C_{\ast}}{\sqrt{n}%
}\right)  \leq2e^{-y}.\label{ExpIneq}%
\end{equation}
As in Comte \cite{Comte}, let $\Theta_{n,b}=\left\{  \left\vert \frac
{\left\Vert \widehat{f}_{n}\right\Vert _{\infty}}{\left\Vert f\right\Vert
_{\infty}}-1\right\vert <b\right\}  $, with $b\in\left(  0,1\right)  $. Then,
using that $P\left(  \left\vert \widehat{b}_{j,k}-b_{j,k}\right\vert
>B+D\right)  \leq$ $P\left(  \left\vert \widehat{b}_{j,k}-b_{j,k}\right\vert
>B\right)  $ for $B,D>0$, and taking $x_{n}=\frac{\delta\log n}{\left(
1-b\right)  ^{2}}$, one gets%
\begin{align*}
&  P\left(  \left\vert \widehat{b}_{j,k}-b_{j,k}\right\vert >\frac
{\widehat{\xi}_{j,n}}{2}\,\right)  \leq P\left(  \left\vert \widehat{b}%
_{j,k}-b_{j,k}\right\vert >2\left\Vert \widehat{f}_{n}\right\Vert _{\infty
}\left(  \sqrt{\frac{x_{n}}{n}}+2^{\frac{j}{2}}\left\Vert \psi\right\Vert
_{\infty}\frac{x_{n}}{n}\right)  \right)  \\
&  \leq P\left(  \left(  \left\vert \widehat{b}_{j,k}-b_{j,k}\right\vert
>2\left\Vert \widehat{f}_{n}\right\Vert _{\infty}\left(  \sqrt{\frac{x_{n}}%
{n}}+2^{\frac{j}{2}}\left\Vert \psi\right\Vert _{\infty}\frac{x_{n}}%
{n}\right)  \right)  \mid\Theta_{n,b}\right)  P\left(  \Theta_{n,b}\right)  \\
&  +P\left(  \left(  \left\vert \widehat{b}_{j,k}-b_{j,k}\right\vert
>2\left\Vert \widehat{f}_{n}\right\Vert _{\infty}\left(  \sqrt{\frac{x_{n}}%
{n}}+2^{\frac{j}{2}}\left\Vert \psi\right\Vert _{\infty}\frac{x_{n}}%
{n}\right)  \right)  \mid\Theta_{n,b}^{c}\right)  P\left(  \Theta_{n,b}%
^{c}\right)  \\
&  :=P_{1}P\left(  \Theta_{n,b}\right)  +P_{2}P\left(  \Theta_{n,b}%
^{c}\right)  .
\end{align*}
In Comte \cite{Comte} is proved that if $\left\Vert f-f_{n}\right\Vert
_{\infty}\leq\frac{1}{4}\left\Vert f\right\Vert _{\infty}$ then $P\left(
\Theta_{n,b}^{c}\right)  \leq O\left(  n^{-4}\right)  $ for the choices of
$1\geq b\approx\frac{4}{6}\sqrt{\frac{5}{\pi}}=0.841\geq\frac{3}{4}$ and
$N_{n}\leq\frac{1}{36\left(  r+1\right)  ^{2}}\frac{n}{\log n}$, where
$\kappa=\frac{1}{36}$ is the numerical constant in the hypothesis of our
theorem. Following its proof it can be shown that this bound can be improved
taking $\kappa=\frac{1}{36\left(  \frac{7}{5}\right)  }$ and $b$ as before
(see the three last equations of page 290 in \cite{Comte}). With this
selection of $\kappa$ we obtain that $P\left(  \Theta_{n,b}^{c}\right)  \leq
O\left(  n^{-6}\right)  $. Using that $P\left(  \Theta_{n,b}\right)  =O\left(
1\right)  $ and $P_{2}=O\left(  1\right)  $, it only remains to bound the
conditional probability $P_{1}$. On $\Theta_{n,b}$ the following inequalities
hold:%
\begin{equation}
\left(  a\right)  \text{ }\left\Vert \widehat{f}_{n}\right\Vert _{\infty
}>\left(  1-b\right)  \left\Vert f\right\Vert _{\infty}\text{ and }\left(
b\right)  \text{ }\left\Vert \widehat{f}_{n}\right\Vert _{\infty}<\left(
1+b\right)  \left\Vert f\right\Vert _{\infty}.\label{IneqI1}%
\end{equation}
Then, using (\ref{IneqI1}$a$) we get
\[
P_{1}\leq P\left(  \left\vert \widehat{b}_{j,k}-b_{j,k}\right\vert
>2\left\Vert f\right\Vert _{\infty}\left(  \sqrt{\frac{\delta\log n}{n}%
}+2^{\frac{j}{2}}\left\Vert \psi\right\Vert _{\infty}\frac{\delta\log n}%
{n}\right)  \right)  \leq2e^{-\delta\log n}=2n^{-\delta},
\]
where the last inequality is obtained using (\ref{ExpIneq}) for $y=\delta\log
n>0$. Hence, using that $\delta=6$, we get $P\left(  \left\vert \widehat
{b}_{j,k}-b_{j,k}\right\vert >\frac{\widehat{\xi}_{j,n}}{2}\right)  \leq
O\left(  n^{-6}\right)  $. Therefore $T^{\prime\prime\prime}\leq
C\sum\limits_{j=j_{0}}^{j_{1}}\sum\limits_{k=0}^{2^{j}-1}n^{-2}\leq O\left(
n^{-2}2^{j_{1}}\right)  \leq O\left(  n^{-\frac{2s}{2s+1}}\right)  .$

Now we bound $T^{\prime}$. Let $j_{A}$ be the integer such that $2^{j_{A}%
}>\left(  \frac{n}{\log n}\right)  ^{\frac{1}{2s+1}}>2^{j_{A}-1}$, then
$T^{\prime}=T_{1}^{\prime}+T_{2}^{\prime}$, where the first component is
computed over the set of indices $j_{0}\leq j\leq j_{A}$ and the second
component over $j_{A}+1\leq j\leq j_{1}$. Hence, using Lemma
\ref{LemmaBoundBiasVarBetahat} we obtain%
\[
T_{1}^{\prime}\leq\sum\limits_{j=j_{0}}^{j_{A}}\sum\limits_{k=0}^{2^{j}%
-1}\mathbb{E}\left(  b_{j,k}-\widehat{b}_{j,k}\right)  ^{2}\leq O\left(
\frac{2^{j_{A}}}{n}\right)  =O\left(  \left(  \frac{n}{\log n}\right)
^{\frac{1}{2s+1}}n^{-1}\right)  \leq O\left(  n^{-\frac{2s}{2s+1}}\right)  .
\]
To bound $T_{2}^{\prime}$, note that $T_{2}^{\prime}:=T_{2,1}^{\prime}%
+T_{2,2}^{\prime}$, where
\begin{align*}
T_{2,1}^{\prime} &  =\sum\limits_{j=j_{A}}^{j_{1}}\sum\limits_{k=0}^{2^{j}%
-1}\mathbb{E}\left\{  \left(  b_{j,k}-\widehat{b}_{j,k}\right)  ^{2}I\left(
\left\vert b_{j,k}\right\vert >\frac{\widehat{\xi}_{j,n}}{2},\Theta
_{n,b}\right)  \right\}  \text{ and}\\
T_{2,2}^{\prime} &  =\sum\limits_{j=j_{A}}^{j_{1}}\sum\limits_{k=0}^{2^{j}%
-1}\mathbb{E}\left\{  \left(  b_{j,k}-\widehat{b}_{j,k}\right)  ^{2}I\left(
\left\vert b_{j,k}\right\vert >\frac{\widehat{\xi}_{j,n}}{2},\Theta_{n,b}%
^{c}\right)  \right\}  .
\end{align*}
Using that on $\Theta_{n,b}$ inequality (\ref{IneqI1}$a$) holds and following
the same procedures as in the proof of Theorem \ref{Teo5.1enAB}, we get the
desired bound for $T_{2,1}^{\prime}$.
\begin{align}
T_{2,1}^{\prime} &  \leq\frac{C}{n}\sum\limits_{j=j_{A}}^{j_{1}}%
\sum\limits_{k=0}^{2^{j}-1}I\left(  \left\vert b_{j,k}\right\vert >2\left(
1-b\right)  \left\Vert f\right\Vert _{\infty}\left(  \sqrt{\frac{x_{n}}{n}%
}+2^{\frac{j}{2}}\left\Vert \psi\right\Vert _{\infty}\frac{x_{n}}{n}\right)
\right)  \label{Line4}\\
&  \leq\frac{C}{4}\frac{1}{n}\sum\limits_{j=j_{A}}^{j_{1}}\sum\limits_{k=0}%
^{2^{j}-1}\frac{\left\vert b_{j,k}\right\vert ^{2}}{\left\Vert f\right\Vert
_{\infty}^{2}\left(  \sqrt{\frac{\delta\log n}{n}}+2^{\frac{j}{2}}\left\Vert
\psi\right\Vert _{\infty}\frac{\delta\log n}{\left(  1-b\right)  n}\right)
^{2}}\nonumber\\
&  \leq O\left(  2^{-2s^{\ast}j_{A}}\right)  =O\left(  \left(  \frac{n}{\log
n}\right)  ^{-\frac{2s^{\ast}}{2s+1}}\right)  ,\nonumber
\end{align}
where we have used that $\sqrt{\delta\log n}+\left\Vert \psi\right\Vert
_{\infty}\left(  1-b\right)  ^{-1}\delta n^{\frac{-s}{2s+1}}\left(  \log
n\right)  ^{\frac{4s+1}{4s+2}}\rightarrow+\infty$ when $n\rightarrow+\infty$
and that condition (\ref{cnd}) is satisfied. Now remark that if $p=2$ then
$s^{\ast}=s$ and thus
\begin{equation}
T_{2,1}^{\prime}=O\left(  \left(  \frac{n}{\log n}\right)  ^{-\frac{2s}{2s+1}%
}\right)  .\label{B1}%
\end{equation}
For the case $1\leq p<2$, from (\ref{Line4}) we have that
\begin{align}
T_{2,1}^{\prime} &  \leq\frac{C}{n}\sum\limits_{j=j_{A}}^{j_{1}}%
\sum\limits_{k=0}^{2^{j}-1}I\left(  \left\vert b_{j,k}\right\vert >2\left\Vert
f\right\Vert _{\infty}\left(  \sqrt{\frac{\delta\log n}{n}}+2^{\frac{j}{2}%
}\left\Vert \psi\right\Vert _{\infty}\frac{\delta\log n}{\left(  1-b\right)
n}\right)  \right)  \nonumber\\
&  \leq\frac{C}{n}\sum\limits_{j=j_{A}}^{j_{1}}\sum\limits_{k=0}^{2^{j}%
-1}\left\vert b_{j,k}\right\vert ^{-p}\left\vert b_{j,k}\right\vert
^{p}I\left(  \left\vert b_{j,k}\right\vert ^{-p}<\left(  2\left\Vert
f\right\Vert _{\infty}\sqrt{\delta\frac{\log n}{n}}\right)  ^{-p}\right)
\nonumber\\
&  \leq\left(  \log n\right)  C\left(  \left\Vert f\right\Vert _{\infty
},\delta,p\right)  \frac{\left(  \log n\right)  ^{-\frac{p}{2}}}{n^{1-\frac
{p}{2}}}\sum\limits_{j=j_{A}}^{j_{1}}\sum\limits_{k=0}^{2^{j}-1}\left\vert
b_{j,k}\right\vert ^{p}=O\left(  \frac{\left(  \log n\right)  ^{1-\frac{p}{2}%
}}{n^{1-\frac{p}{2}}}2^{-pj_{A}s^{\ast}}\right)  \nonumber\\
&  \leq O\left(  \left(  \frac{n}{\log n}\right)  ^{-\frac{2s}{2s+1}}\right)
,\label{B2sust}%
\end{align}
where we have used condition (\ref{cnd1}). Hence $T_{2,1}^{\prime}=O\left(
\left(  \frac{n}{\log n}\right)  ^{-\frac{2s}{2s+1}}\right)  $.

Now we bound $T_{2,2}^{\prime}$. Using Cauchy-Schwarz inequality, we have
\begin{align}
T_{2,2}^{\prime}  &  \leq C\sum\limits_{j=j_{A}}^{j_{1}}\sum\limits_{k=0}%
^{2^{j}-1}nP^{\frac{1}{2}}\left(  \left\vert b_{j,k}\right\vert >2\left\Vert
\widehat{f}_{n}\right\Vert _{\infty}\left(  \sqrt{\frac{x_{n}}{n}}+2^{\frac
{j}{2}}\left\Vert \psi\right\Vert _{\infty}\frac{x_{n}}{n}\right)
+\sqrt{\frac{\log n}{n}}\mid\Theta_{n,b}^{c}\right)  P^{\frac{1}{2}}\left(
\Theta_{n,b}^{c}\right) \nonumber\\
&  \leq C\sum\limits_{j=j_{A}}^{j_{1}}\sum\limits_{k=0}^{2^{j}-1}nP^{\frac
{1}{2}}\left(  \Theta_{n,b}^{c}\right)  \leq C\sum\limits_{j=j_{A}}^{j_{1}%
}\sum\limits_{k=0}^{2^{j}-1}n^{-2}\leq O\left(  \frac{2^{j_{1}}}{n^{2}%
}\right)  \leq O\left(  n^{-\frac{2s}{2s+1}}\right)  , \label{B3}%
\end{align}
where we have used that $\mathbb{E}\left\{  \left(  b_{j,k}-\widehat{b}%
_{j,k}\right)  ^{4}\right\}  =O\left(  n^{2}\right)  $ and that $P\left(
\Theta_{n,b}^{c}\right)  \leq O\left(  n^{-6}\right)  $. Then, putting
together (\ref{B1}), (\ref{B2sust}) and (\ref{B3}), we obtain that
$T_{2}^{\prime}=O\left(  \left(  \frac{n}{\log n}\right)  ^{-\frac{2s}{2s+1}%
}\right)  $.

Now we bound $T^{\prime\prime}$. Set $j_{A}$ as before, then $T^{\prime\prime
}$ $=T_{1}^{\prime\prime}+T_{2}^{\prime\prime}$, where the first component is
calculated over the set of indices $j_{0}\leq j\leq j_{A}$ and the second
component over $j_{A}+1\leq j\leq j_{1}$. Recall that $x_{n}=\frac{\delta\log
n}{\left(  1-b\right)  ^{2}}$, then $T_{1}^{\prime\prime}\leq T_{1,1}%
^{\prime\prime}+T_{1,2}^{\prime\prime}$, where
\begin{align*}
T_{1,1}^{\prime\prime} &  =\sum\limits_{j=j_{0}}^{j_{A}}\sum\limits_{k=0}%
^{2^{j}-1}b_{j,k}^{2}P\left(  \left\vert b_{j,k}\right\vert \leq4\left[
2\left(  1+b\right)  \left\Vert f\right\Vert _{\infty}\left(  \sqrt
{\frac{x_{n}}{n}}+2^{\frac{j}{2}}\left\Vert \psi\right\Vert _{\infty}%
\frac{x_{n}}{n}\right)  +\sqrt{\frac{\log n}{n}}\right]  \right)  \text{
and}\\
T_{1,2}^{\prime\prime} &  =\sum\limits_{j=j_{0}}^{j_{A}}\sum\limits_{k=0}%
^{2^{j}-1}b_{j,k}^{2}P\left(  \Theta_{n,b}^{c}\right)  ,
\end{align*}
where we have used that given $\Theta_{n,b}$ inequality (\ref{IneqI1}$b$)
holds. For $T_{1,1}^{\prime\prime}$ we have
\begin{align}
T_{1,1}^{\prime\prime} &  =\sum\limits_{j=j_{0}}^{j_{A}}\sum\limits_{k=0}%
^{2^{j}-1}b_{j,k}^{2}P\left(  \left\vert b_{j,k}\right\vert ^{2}\leq16\left[
2\left(  1+b\right)  \left\Vert f\right\Vert _{\infty}\left(  \sqrt
{\frac{x_{n}}{n}}+2^{\frac{j}{2}}\left\Vert \psi\right\Vert _{\infty}%
\frac{x_{n}}{n}\right)  +\sqrt{\frac{\log n}{n}}\right]  ^{2}\right)
\nonumber\\
&  \leq16\sum\limits_{j=j_{0}}^{j_{A}}\sum\limits_{k=0}^{2^{j}-1}\left[
2\left(  1+b\right)  \left\Vert f\right\Vert _{\infty}\left(  \sqrt
{\frac{\delta\log n}{\left(  1-b\right)  ^{2}n}}+2^{\frac{j}{2}}\left\Vert
\psi\right\Vert _{\infty}\frac{\delta\log n}{\left(  1-b\right)  ^{2}%
n}\right)  +\sqrt{\frac{\log n}{n}}\right]  ^{2}\nonumber\\
&  \leq C\sum\limits_{j=j_{0}}^{j_{A}}\sum\limits_{k=0}^{2^{j}-1}\left(
C\left(  \left\Vert f\right\Vert _{\infty},b\right)  \left(  \frac{\delta\log
n}{\left(  1-b\right)  ^{2}n}+\left\Vert \psi\right\Vert _{\infty}^{2}%
\frac{\delta^{2}\log n}{\left(  1-b\right)  ^{4}n}\right)  +\frac{\log n}%
{n}\right)  \nonumber\\
&  =O\left(  \left(  \frac{n}{\log n}\right)  ^{-\frac{2s}{2s+1}}\right)
,\label{B5}%
\end{align}
where we have used repeatedly that $\left(  B+D\right)  ^{2}\leq2\left(
B^{2}+D^{2}\right)  $ for all $B,D\in\mathbb{R}$.

To bound $T_{1,2}^{\prime\prime}$ we use again that $P\left(  \Theta_{n,b}%
^{c}\right)  \leq O\left(  n^{-6}\right)  $ and that condition (\ref{cnd}) is
satisfied. Then
\begin{equation}
T_{1,2}^{\prime\prime}\leq\sum\limits_{j=j_{0}}^{j_{A}}\sum\limits_{k=0}%
^{2^{j}-1}b_{j,k}^{2}n^{-6}\leq n^{-6}\sum\limits_{j=j_{0}}^{j_{A}%
}C2^{-2js^{\ast}}=O\left(  n^{-6}2^{-2j_{0}s^{\ast}}\right)  \leq O\left(
n^{-1}\right)  .\label{B6}%
\end{equation}
Hence, by (\ref{B5}) and (\ref{B6}), $T_{1}^{\prime\prime}=O\left(  \left(
\frac{n}{\log n}\right)  ^{-\frac{2s}{2s+1}}\right)  $. Now we bound
$T_{2}^{\prime\prime}$.
\[
T_{2}^{\prime\prime}\leq\sum\limits_{j=j_{A}+1}^{j_{1}}\sum\limits_{k=0}%
^{2^{j}-1}b_{j,k}^{2}P\left(  \left\vert b_{j,k}\right\vert \leq2\widehat{\xi
}_{j,n}\,\right)  \leq\sum\limits_{j=j_{A}+1}^{j_{1}}\sum\limits_{k=0}%
^{2^{j}-1}b_{j,k}^{2}=O\left(  2^{-2j_{A}s^{\ast}}\right)  =O\left(  \left(
\frac{n}{\log n}\right)  ^{-\frac{2s^{\ast}}{2s+1}}\right)  ,
\]
where we have used again the condition (\ref{cnd}). Now remark that if $p=2$
then $s^{\ast}=s$ and thus $T_{2}^{\prime\prime}=O\left(  \left(  \frac
{n}{\log n}\right)  ^{-\frac{2s^{\ast}}{2s+1}}\right)  =O\left(  \left(
\frac{n}{\log n}\right)  ^{-\frac{2s}{2s+1}}\right)  $. For $1\leq p<2$, we
proceed as follows.%
\begin{align*}
T_{2}^{\prime\prime} &  =\sum\limits_{j=j_{A}+1}^{j_{1}}\sum\limits_{k=0}%
^{2^{j}-1}\mathbb{E}\left[  b_{j,k}^{2}I\left(  \left\vert b_{j,k}\right\vert
\leq2\widehat{\xi}_{j,n}\,,\Theta_{n,b}\right)  +b_{j,k}^{2}I\left(
\left\vert b_{j,k}\right\vert \leq2\widehat{\xi}_{j,n}\,,\Theta_{n,b}%
^{c}\right)  \right]  \\
&  \leq\sum\limits_{j=j_{A}+1}^{j_{1}}\sum\limits_{k=0}^{2^{j}-1}%
\mathbb{E}\left[  b_{j,k}^{2}I\left(  \left\vert b_{j,k}\right\vert
\leq4\left(  2\left(  1+b\right)  \left\Vert f\right\Vert _{\infty}\left(
\sqrt{\frac{x_{n}}{n}}+2^{\frac{j}{2}}\left\Vert \psi\right\Vert _{\infty
}\frac{x_{n}}{n}\right)  +\sqrt{\frac{\log n}{n}}\right)  \right)  \right]  \\
&  +\sum\limits_{j=j_{A}+1}^{j_{1}}\sum\limits_{k=0}^{2^{j}-1}b_{j,k}%
^{2}P\left(  \left\vert b_{j,k}\right\vert \leq\left(  8\left\Vert \widehat
{f}_{n}\right\Vert _{\infty}\left(  \sqrt{\frac{x_{n}}{n}}+2^{\frac{j}{2}%
}\left\Vert \psi\right\Vert _{\infty}\frac{x_{n}}{n}\right)  +4\sqrt
{\frac{\log n}{n}}\right)  \mid\Theta_{n,b}^{c}\right)  P\left(  \Theta
_{n,b}^{c}\right)  \\
&  :=T_{2,1}^{\prime\prime}+T_{2,2}^{\prime\prime},
\end{align*}
where we have used that on $\Theta_{n,b}$ inequality (\ref{IneqI1}$b$) holds.
Now we bound $T_{2,1}^{\prime\prime}$.
\begin{align}
T_{2,1}^{\prime\prime} &  \leq\sum\limits_{j=j_{A}+1}^{j_{1}}\sum
\limits_{k=0}^{2^{j}-1}\left\vert b_{j,k}\right\vert ^{2-p}\left\vert
b_{j,k}\right\vert ^{p}I\left(  \left\vert b_{j,k}\right\vert \leq8\left(
1+b\right)  \left\Vert f\right\Vert _{\infty}\left(  \sqrt{\frac{x_{n}}{n}%
}+2^{\frac{j}{2}}\left\Vert \psi\right\Vert _{\infty}\frac{x_{n}}{n}\right)
+4\sqrt{\frac{\log n}{n}}\right)  \nonumber\\
&  \leq\sum\limits_{j=j_{A}+1}^{j_{1}}\sum\limits_{k=0}^{2^{j}-1}\left(
8\left(  1+b\right)  \left\Vert f\right\Vert _{\infty}\left(  \sqrt
{\frac{\delta\log n}{\left(  1-b\right)  ^{2}n}}+2^{\frac{j_{1}}{2}}%
\frac{\left\Vert \psi\right\Vert _{\infty}\delta\log n}{\left(  1-b\right)
^{2}n}\right)  +4\sqrt{\frac{\log n}{n}}\right)  ^{2-p}\left\vert
b_{j,k}\right\vert ^{p}\nonumber\\
&  \leq C\left(  \left\Vert f\right\Vert _{\infty},b,\delta,\left\Vert
\psi\right\Vert _{\infty}\right)  ^{2-p}\left(  \sqrt{\frac{\log n}{n}%
}\right)  ^{2-p}\sum\limits_{j=j_{A}+1}^{j_{1}}\sum\limits_{k=0}^{2^{j}%
-1}\left\vert b_{j,k}\right\vert ^{p}\leq O\left(  \left(  \frac{\log n}%
{n}\right)  ^{\frac{2-p}{2}}2^{-pj_{A}s^{\ast}}\right)  \nonumber\\
&  =O\left(  \left(  \frac{\log n}{n}\right)  ^{\frac{2-p}{2}}\left(  \frac
{n}{\log n}\right)  ^{-\frac{p\left(  s+\frac{1}{2}-\frac{1}{p}\right)
}{2s+1}}\right)  =O\left(  \left(  \frac{n}{\log n}\right)  ^{-\frac{2s}%
{2s+1}}\right)  .\label{B8}%
\end{align}
where we have used that condition (\ref{cnd1}) is satisfied. To bound
$T_{2,2}^{\prime\prime}$ we use again that $P\left(  \Theta_{n,b}^{c}\right)
\leq O\left(  n^{-6}\right)  $ and that condition (\ref{cnd}) also holds.
Then, from (\ref{B6}) we get%
\begin{equation}
T_{2,2}^{\prime\prime}\leq\sum\limits_{j=j_{A}+1}^{j_{1}}\sum\limits_{k=0}%
^{2^{j}-1}b_{j,k}^{2}P\left(  \Theta_{n,b}^{c}\right)  \leq n^{-6}%
\sum\limits_{j=j_{A}+1}^{j_{1}}C2^{-2js^{\ast}}=O\left(  n^{-6}2^{-2j_{A}%
s^{\ast}}\right)  \leq O\left(  n^{-1}\right)  .\label{B7}%
\end{equation}
Hence, by (\ref{B8}) and (\ref{B7}), $T^{\prime\prime}=O\left(  \left(
\frac{n}{\log n}\right)  ^{-\frac{2s}{2s+1}}\right)  $. Combining all terms in
(\ref{Sum4terminos}), we conclude that:%
\[
\mathbb{E}\left\Vert \beta-\widehat{\beta}_{\widehat{\xi}_{j,n}}\right\Vert
_{2}^{2}=O\left(  \left(  \frac{n}{\log n}\right)  ^{-\frac{2s}{2s+1}}\right)
.
\]
This completes the proof.
\end{proof}

\subsection{Proof of Theorem \ref{theo:lin}}

First, one needs the following proposition.

\begin{proposition}
\label{prop:existlin} \label{ExistenceEstimator} Let $\beta_{j,k}=\left\langle
f,\psi_{j,k}\right\rangle $ and $\widehat{\beta}_{j,k}=\left\langle I_{n}%
,\psi_{j,k}\right\rangle $ with $\left(  j,k\right)  \in$ $\Lambda_{j_{1}}$.
Suppose that $f\in F_{p,q}^{s}\left(  M\right)  $ with $s>1/p$ and $1\leq
p\leq2$. Let $M_{1}>0$ be a constant such that $M_{1}^{-1}\leq$ $f\leq$
$M_{1}$ (see Lemma \ref{LemaA.4inAB}). Let $\epsilon_{j_{1}}=2M_{1}%
^{2}e^{2\gamma_{j_{1}}+1}D_{j_{1}}A_{j_{1}}$. If $\epsilon_{j_{1}}\leq1$, then
there exists $\theta_{j_{1}}^{\ast}\in\mathbb{R}^{\#\Lambda_{j_{1}}}$ such
that:%
\[
\left\langle f_{j_{1},\theta_{j_{1}}^{\ast}},\psi_{j,k}\right\rangle
=\left\langle f,\psi_{j,k}\right\rangle =\beta_{j,k}\text{ for all }\left(
j,k\right)  \in\Lambda_{j_{1}}%
\]
Moreover, the following inequality holds (approximation error)
\[
\Delta\left(  f;f_{j_{1},\theta_{j_{1}}^{\ast}}\right)  \leq\frac{M_{1}}%
{2}e^{\gamma_{j_{1}}}D_{j_{1}}^{2}.
\]
Suppose that Assumptions \ref{ass1} and \ref{ass2} hold. Let $\eta_{j_{1}%
,n}=4M_{1}^{2}e^{2\gamma_{j_{1}}+2\epsilon_{j_{1}}+2}A_{j_{1}}^{2}%
\frac{\#\Lambda_{j_{1}}}{n}$. Then, for every $\lambda>0$ such that
$\lambda\leq\eta_{j_{1},n}^{-1}$ there exists a set $\Omega_{n,1}$ of
probability less than $M_{2}\lambda^{-1}$, where $M_{2}$ is the constant
defined in Lemma \ref{LemmaBoundBiasVarBetahat}, such that outside the set
$\Omega_{n,1}$ there exists some $\widehat{\theta}_{n}\in\mathbb{R}%
^{\#\Lambda_{j_{1}}}$ which satisfies:%
\[
\left\langle f_{j_{1},\widehat{\theta}_{n}},\psi_{j,k}\right\rangle
=\left\langle I_{n},\psi_{j,k}\right\rangle =\widehat{\beta}_{j,k}\text{ for
all }\left(  j,k\right)  \in\Lambda_{j_{1}}.
\]
Moreover, outside the set $\Omega_{n,1}$, the following inequality holds
(estimation error)
\[
\Delta\left(  f_{j_{1},\theta_{j_{1}}^{\ast}};f_{j_{1},\widehat{\theta}_{n}%
}\right)  \leq2M_{1}e^{\gamma_{j_{1}}+\epsilon_{j_{1}}+1}M_{2}\lambda
\frac{\#\Lambda_{j_{1}}}{n}.
\]

\end{proposition}

\begin{proof}
\noindent\textbf{Approximation error:} Recall that $\beta_{j,k}=\left\langle
f,\psi_{j,k}\right\rangle $ and let $\beta=\left(  \beta_{j,k}\right)
_{(j,k)\in\Lambda_{j_{1}}}$. Define by $g_{j_{1}}=\sum\limits_{(j,k)\in
\Lambda_{j_{1}}}\theta_{j,k}\psi_{j,k}$ an approximation of $g=\log\left(
f\right)  $ and let $\beta_{0,\left(  j,k\right)  }=\left\langle
f_{j_{1},\theta_{j_{1}}},\psi_{j,k}\right\rangle =\left\langle \exp\left(
g_{j_{1}}\right)  ,\psi_{j,k}\right\rangle $ with $\theta_{j_{1}}=\left(
\theta_{j,k}\right)  _{(j,k)\in\Lambda_{j_{1}}}$ and $\beta_{0}=\left(
\beta_{0,\left(  j,k\right)  }\right)  _{(j,k)\in\Lambda_{j_{1}}}$. Observe
that the coefficients $\beta_{j,k}-\beta_{0,\left(  j,k\right)  }$,
$(j,k)\in\Lambda_{j_{1}}$, are the coefficients of the orthonormal projection
of $f-f_{j_{1},\theta_{j_{1}}}$ onto $V_{j}$. Hence by Bessel's inequality,
$\left\Vert \beta-\beta_{0}\right\Vert _{2}^{2}\leq\left\Vert f-f_{j_{1}%
,\theta_{j_{1}}}\right\Vert _{L_{2}}^{2}$. Using Lemma \ref{LemaA.4inAB} and
Lemma 2 in Barron and Sheu \cite{B-S}, we get that:%
\begin{align*}
\left\Vert \beta-\beta_{0}\right\Vert _{2}^{2} &  \leq\int\left(
f-f_{j_{1},\theta_{j_{1}}}\right)  ^{2}d\mu\leq M_{1}\int\frac{\left(
f-f_{j_{1},\theta_{j_{1}}}\right)  ^{2}}{f}d\mu\\
&  \leq M_{1}e^{2\left\Vert \log\left(  \frac{f}{f_{j_{1},\theta_{j_{1}}}%
}\right)  \right\Vert _{\infty}}\int f\left(  \log\left(  \frac{f}%
{f_{j_{1},\theta_{j_{1}}}}\right)  \right)  ^{2}d\mu\\
&  \leq M_{1}^{2}e^{2\left\Vert g-g_{j_{1}}\right\Vert _{\infty}}\left\Vert
g-g_{j_{1}}\right\Vert _{L_{2}}^{2}=M_{1}^{2}e^{2\gamma_{j_{1}}}D_{j_{1}}^{2}.
\end{align*}
Then, one can easily check that $b=e^{\left(  \left\Vert \log\left(
f_{j_{1},\theta_{j_{1}}}\right)  \right\Vert _{\infty}\right)  }\leq
M_{1}e^{\gamma_{j_{1}}}$. Thus the assumption that $\epsilon_{j_{1}}\leq1$
implies that the inequality $\left\Vert \beta-\beta_{0}\right\Vert _{2}\leq
M_{1}e^{\gamma_{j_{1}}}D_{j_{1}}\leq\frac{1}{2beA_{j_{1}}}$ is satisfied.
Hence, Lemma \ref{Lema3.2enAB} can be applied with $\theta_{0}=\theta_{j_{1}}%
$, $\widetilde{\beta}=\beta$ and $b=\exp\left(  \left\Vert \log\left(
f_{j_{1},\theta_{j_{1}}}\right)  \right\Vert _{\infty}\right)  $, which
implies that there exists $\theta_{j_{1}}^{\ast}=\theta\left(  \beta\right)  $
such that $\left\langle f_{j_{1},\theta_{j_{1}}^{\ast}},\psi_{j,k}%
\right\rangle =\beta_{j,k}$ for all $\left(  j,k\right)  \in$ $\Lambda_{j_{1}%
}$.

By the Pythagorian-like relationship (\ref{IdPyth}), we obtain that
$\Delta\left(  f;f_{j_{1},\theta_{j_{1}}^{\ast}}\right)  \leq\Delta\left(
f;f_{j_{1},\theta_{j_{1}}}\right)  .$ Now we use a result wich states that if
$f$ and $g$ are two functions in $L_{2}([0,1])$ such that $\log\left(
\frac{f}{g}\right)  $ is bounded. Then $\Delta\left(  f;g\right)  \leq\frac
{1}{2}e^{\left\Vert \log\left(  \frac{f}{g}\right)  \right\Vert _{\infty}}%
\int_{0}^{1}f\left(  \log\left(  \frac{f}{g}\right)  \right)  ^{2}d\mu$, where
$\mu$ denotes the Lebesgue measure on $[0,1]$.(see Lemma A.1 in Antoniadis and
Bigot \cite{A-B}). Hence, it follows that
\begin{align*}
\Delta\left(  f;f_{j_{1},\theta_{j_{1}}^{\ast}}\right)   &  \leq\frac{1}%
{2}e^{\left\Vert \log\left(  \frac{f}{f_{j_{1},\theta_{j_{1}}}}\right)
\right\Vert _{\infty}}\int f\left(  \log\left(  \frac{f}{f_{j_{1}%
,\theta_{j_{1}}}}\right)  \right)  ^{2}d\mu\\
&  =\frac{M_{1}}{2}e^{\left\Vert g-g_{j_{1}}\right\Vert _{\infty}}\left\Vert
g-g_{j_{1}}\right\Vert _{L_{2}}^{2}=\frac{M_{1}}{2}e^{\gamma_{j_{1}}}D_{j_{1}%
}^{2}.
\end{align*}
which completes the proof for the approximation error. \newline


\noindent\textbf{Estimation error:} Applying again Lemma \ref{Lema3.2enAB}
with $\theta_{0}=\theta_{j_{1}}^{\ast}$, $\beta_{0,\left(  j,k\right)  }=$
$\left\langle f_{j_{1},\theta_{0}},\psi_{j,k}\right\rangle =\beta_{j,k}$,
$\widetilde{\beta}=\widehat{\beta}$, where $\widehat{\beta}=\left(
\widehat{\beta}_{j,k}\right)  _{\left(  j,k\right)  \in\Lambda_{j_{1}}}$, and
$b=\exp\left(  \left\Vert \log\left(  f_{j_{1},\theta_{j_{1}}^{\ast}}\right)
\right\Vert _{\infty}\right)  $ we obtain that if $\left\Vert \widehat{\beta
}-\beta\right\Vert _{2}\leq\frac{1}{2ebA_{j_{1}}}$ with $\beta=\left(
\beta_{j,k}\right)  _{\left(  j,k\right)  \in\Lambda_{j_{1}}}$ then there
exists $\widehat{\theta}_{n}=\theta\left(  \widehat{\beta}\right)  $ such that
$\left\langle f_{j_{1},\widehat{\theta}_{n}},\psi_{j,k}\right\rangle
=\widehat{\beta}_{j,k}$ for all $\left(  j,k\right)  \in$ $\Lambda_{j_{1}}$.

Hence, it remains to prove that our assumptions imply that the event
$\left\Vert \widehat{\beta}-\beta\right\Vert _{2}\leq\frac{1}{2ebA_{j_{1}}}$
holds with probability $1-M_{2}\lambda^{-1}$. First remark that $b\leq
M_{1}e^{\gamma_{j_{1}}+\epsilon_{j_{1}}}$ and that by Markov's inequality and
Lemma \ref{LemmaBoundBiasVarBetahat} we obtain that for any $\lambda>0$,
$P\left(  \left\Vert \widehat{\beta}-\beta\right\Vert _{2}^{2}\geq\lambda
\frac{\#\Lambda_{j_{1}}}{n}\right)  \leq\frac{1}{\lambda}\frac{n}%
{\#\Lambda_{j_{1}}}\mathbb{E}\left\Vert \widehat{\beta}-\beta\right\Vert
_{2}^{2}\leq M_{2}\lambda^{-1}$. Hence, outside a set $\Omega_{n,1}$ of
probability less than $M_{2}\lambda^{-1}$ then $\left\Vert \widehat{\beta
}-\beta\right\Vert _{2}^{2}\leq\lambda\frac{\#\Lambda_{j_{1}}}{n}.$ Therefore,
the condition $\left\Vert \widehat{\beta}-\beta\right\Vert _{2}\leq\frac
{1}{2ebA_{j_{1}}}$ holds if $\left(  \lambda\frac{\#\Lambda_{j_{1}}}%
{n}\right)  ^{\frac{1}{2}}\leq\frac{1}{2ebA_{j_{1}}}$, which is equivalent to
$4e^{2}b^{2}A_{j_{1}}^{2}\lambda\frac{\#\Lambda_{j_{1}}}{n}\leq1$. This last
inequality is true if $\eta_{j_{1},n}=4M_{1}^{2}e^{2\gamma_{j_{1}}%
+2\epsilon_{j_{1}}+2}A_{j_{1}}^{2}\frac{\#\Lambda_{j_{1}}}{n}\leq\frac
{1}{\lambda}$, using that $b^{2}\leq M_{1}^{2}e^{2\gamma_{j_{1}}%
+2\epsilon_{j_{1}}}$.

Hence, outside the set $\Omega_{n,1}$, our assumptions imply that there exists
$\widehat{\theta}_{n}=\theta\left(  \widehat{\beta}\right)  $ such that
$\left\langle f_{j_{1},\widehat{\theta}_{n}},\psi_{j,k}\right\rangle
=\widehat{\beta}_{j,k}$ for all $\left(  j,k\right)  \in$ $\Lambda_{j_{1}}$.
Finally, outside the set $\Omega_{n,1}$, by using the bound given in Lemma
\ref{Lema3.2enAB}, one obtains the following inequality for the estimation
error%
\[
\Delta\left(  f_{j_{1},\theta_{j_{1}}^{\ast}};f_{j_{1},\widehat{\theta}_{n}%
}\right)  \leq2M_{1}e^{\gamma_{j_{1}}+\epsilon_{j_{1}}+1}\lambda
\frac{\#\Lambda_{j_{1}}}{n}.
\]
which completes the proof of Proposition \ref{prop:existlin}.
\end{proof}

Our assumptions on $j_{1}(n)$ imply that $\frac{1}{2}n^{\frac{1}{2s+1}}%
\leq2^{j_{1}(n)}\leq n^{\frac{1}{2s+1}}$. Therefore, using Lemma
\ref{LemmaAjDj}, one has that for all $f\in F_{2,2}^{s}(M)$ with $s>1/2$
\[
\gamma_{j_{1}\left(  n\right)  }\leq Cn^{\frac{1-2s}{2\left(  2s+1\right)  }%
},\quad A_{j_{1}\left(  n\right)  }\leq Cn^{\frac{1}{2(2s+1)}},\quad
D_{j_{1}\left(  n\right)  }\leq Cn^{-\frac{s}{2s+1}},
\]
where $C$ denotes constants not depending on $g=\log(f)$. Hence,
$\underset{n\rightarrow+\infty}{\lim}\epsilon_{j_{1}\left(  n\right)
}=\underset{n\rightarrow+\infty}{\lim}2M_{1}^{2}e^{2\gamma_{j_{1}\left(
n\right)  }+1}A_{j_{1}\left(  n\right)  }D_{j_{1}\left(  n\right)  }%
=0,$uniformly over $F_{2,2}^{s}(M)$ for $s>1/2$. For all sufficiently large
$n$, $\epsilon_{j_{1}\left(  n\right)  }\leq1$ and thus, using Proposition
\ref{prop:existlin}, there exists $\theta_{j_{1}(n)}^{\ast}\in{\mathbb{R}%
}^{\#\Lambda_{j_{1}(n)}}$ such that
\begin{equation}
\Delta\left(  f;f_{j,\theta_{j_{1}(n)}^{\ast}}\right)  \leq\frac{M_{1}}%
{2}e^{\gamma_{j_{1}(n)}}D_{j_{1}(n)}^{2}\leq Cn^{-\frac{2s}{2s+1}%
}\mbox{ for all }f\in F_{2,2}^{s}(M).\label{eq:deltalin1}%
\end{equation}
By the same arguments it follows that $\underset{n\rightarrow+\infty}{\lim
}\eta_{j_{1}(n),n}=\underset{n\rightarrow+\infty}{\lim}4M_{1}^{2}%
e^{2\gamma_{j_{1}(n)}+2\epsilon_{j_{1}(n)}+2}A_{j_{1}(n)}^{2}\frac
{\#\Lambda_{j_{1}(n)}}{n}=0$, uniformly over $F_{2,2}^{s}(M)$ for $s>1/2$. Now
let $\lambda>0$. The above result shows that for sufficiently large $n$,
$\lambda\leq\eta_{j_{1}(n),n}^{-1}$, and thus using Proposition
\ref{prop:existlin} it follows that there exists a set $\Omega_{n,1}$ of
probability less than $M_{2}\lambda^{-1}$ such that outside this set there
exists $\widehat{\theta}_{n}\in\mathbb{R}^{\#\Lambda_{j_{1}(n)}}$ which
satisfies:
\begin{equation}
\Delta\left(  f_{j_{1}(n),\theta_{j_{1}(n)}^{\ast}};f_{j_{1}(n),\widehat
{\theta}_{n}}\right)  \leq2M_{1}e^{\gamma_{j_{1}(n)}+\epsilon_{j_{1}(n)}%
+1}M_{2}\lambda\frac{\#\Lambda_{j_{1}(n)}}{n}\leq C\lambda n^{-\frac{2s}%
{2s+1}},\label{eq:deltalin2}%
\end{equation}
for all $f\in F_{2,2}^{s}(M)$. Then, by the Pythagorian-like identity
(\ref{IdPyth}) it follows that outside the set $\Omega_{n,1}$
\[
\Delta\left(  f;f_{j_{1}(n),\widehat{\theta}_{n}}\right)  =\Delta\left(
f;f_{j_{1}(n),\theta_{j_{1}(n)}^{\ast}}\right)  +\Delta\left(  f_{j_{1}%
(n),\theta_{j_{1}(n)}^{\ast}};f_{j_{1}(n),\widehat{\theta}_{n}}\right)  ,
\]
and thus Theorem \ref{theo:lin} follows from inequalities (\ref{eq:deltalin1})
and (\ref{eq:deltalin2}).

\subsection{Proof of Theorem \ref{Teo5.1enAB}}

First, one needs the following proposition.

\begin{proposition}
\label{prop:existnonlin} \label{ExistenceAdaptEstimator} Let $\beta
_{j,k}:=\left\langle f,\psi_{j,k}\right\rangle $ and $\widehat{\beta}%
_{\xi_{j,n},\left(  j,k\right)  }:=\delta_{\xi_{j,n}}\left(  \widehat{\beta
}_{j,k}\right)  $ with $\left(  j,k\right)  \in$ $\Lambda_{j_{1}}$. Assume
that $f\in F_{p,q}^{s}\left(  A\right)  $ with $s>1/p$ and $1\leq p\leq2$. Let
$M_{1}>0$ be a constant such that $M_{1}^{-1}\leq$ $f\leq$ $M_{1}$ (see Lemma
\ref{LemaA.4inAB}). Let $\epsilon_{j_{1}}=2M_{1}^{2}e^{2\gamma_{j_{1}}%
+1}D_{j_{1}}A_{j_{1}}$. If $\epsilon_{j_{1}}\leq1$, then there exists
$\theta_{j_{1}}^{\ast}\in\mathbb{R}^{\#\Lambda_{j_{1}}}$ such that:%
\[
\left\langle f_{j_{1},\theta_{j_{1}}^{\ast}},\psi_{j,k}\right\rangle
=\left\langle f,\psi_{j,k}\right\rangle =\beta_{j,k}\text{ for all }\left(
j,k\right)  \in\Lambda_{j_{1}}%
\]
Moreover, the following inequality holds (approximation error)
\[
\Delta\left(  f;f_{j_{1},\theta_{j_{1}}^{\ast}}\right)  \leq\frac{M_{1}}%
{2}e^{\gamma_{j_{1}}}D_{j_{1}}^{2}.
\]
Suppose that Assumptions \ref{ass1} and \ref{ass2} hold. Let $\eta_{j_{1}%
,n}=4M_{1}^{2}e^{2\gamma_{j_{1}}+2\epsilon_{j_{1}}+2}A_{j_{1}}^{2}\left(
\frac{n}{\log n}\right)  ^{-\frac{2s}{2s+1}}$. Then, for every $\lambda>0$
such that $\lambda\leq\eta_{j_{1},n}^{-1}$ there exists a set $\Omega_{n,2}$
of probability less than $M_{3}\lambda^{-1}$, where $M_{3}$ is the constant
defined in Lemma \ref{LemmaBoundBiasVarThresBetahat}, such that outside the
set $\Omega_{n,2}$ there exists some $\widehat{\theta}_{n}\in\mathbb{R}%
^{\#\Lambda_{j_{1}}}$ which satisfies:%
\[
\left\langle f_{j_{1},\widehat{\theta}_{n},\xi_{j,n}}^{HT},\psi_{j,k}%
\right\rangle =\delta_{\xi_{j,n}}\left(  \widehat{\beta}_{j,k}\right)
=\widehat{\beta}_{\xi_{j,n},\left(  j,k\right)  }\text{ for all }\left(
j,k\right)  \in\Lambda_{j_{1}}.
\]
Moreover, outside the set $\Omega_{n,2}$, the following inequality holds
(estimation error)
\[
\Delta\left(  f_{j_{1},\theta_{j_{1}}^{\ast}};f_{j_{1},\widehat{\theta}%
_{n},\xi_{j,n}}^{HT}\right)  \leq2M_{1}e^{\gamma_{j_{1}}+\epsilon_{j_{1}}%
+1}\lambda\left(  \frac{n}{\log n}\right)  ^{-\frac{2s}{2s+1}}.
\]

\end{proposition}

\begin{proof}
\noindent\textbf{Approximation error:} The proof is the same that the one of
Proposition \ref{prop:existlin}.


\noindent\textbf{Estimation error:} Applying Lemma \ref{Lema3.2enAB} with
$\theta_{0}=\theta_{j_{1}}^{\ast}$, $\beta_{0,\left(  j,k\right)  }=$
$\left\langle f_{j_{1},\theta_{0}},\psi_{j,k}\right\rangle =\beta_{j,k}$,
$\widetilde{\beta}=\widehat{\beta}_{\xi_{j,n}}$, where $\widehat{\beta}%
_{\xi_{j,n}}=\left(  \widehat{\beta}_{\xi_{j,n},\left(  j,k\right)  }\right)
_{\left(  j,k\right)  \in\Lambda_{j_{1}}}$, and $b=\exp\left(  \left\Vert
\log\left(  f_{j_{1},\theta_{j_{1}}^{\ast}}\right)  \right\Vert _{\infty
}\right)  $ we obtain that if $\left\Vert \widehat{\beta}_{\xi_{j,n}}%
-\beta\right\Vert _{2}\leq\frac{1}{2ebA_{j_{1}}}$ with $\beta=\left(
\beta_{j,k}\right)  _{\left(  j,k\right)  \in\Lambda_{j_{1}}}$then there
exists $\widehat{\theta}_{n}=\theta\left(  \widehat{\beta}_{\xi_{j,n}}\right)
$ such that $\left\langle f_{j_{1},\widehat{\theta}_{j_{1}},\xi_{j,n}}%
^{HT},\psi_{j,k}\right\rangle =\widehat{\beta}_{\xi_{j,n},\left(  j,k\right)
}$ for all $\left(  j,k\right)  \in$ $\Lambda_{j_{1}}$.

Hence, it remains to prove that our assumptions imply that the event
$\left\Vert \widehat{\beta}_{\xi_{j,n}}-\beta\right\Vert _{2}\leq\frac
{1}{2ebA_{j_{1}}}$ holds with probability $1-M_{3}\lambda^{-1}$. First remark
that $b\leq M_{1}e^{\gamma_{j_{1}}+\epsilon_{j_{1}}}$ and that by Markov's
inequality and Lemma \ref{LemmaBoundBiasVarThresBetahat} we obtain that for
any $\lambda>0$,
\begin{align*}
P\left(  \left\Vert \widehat{\beta}_{\xi_{j,n}}-\beta\right\Vert _{2}^{2}%
\geq\lambda\left(  \frac{n}{\log n}\right)  ^{-\frac{2s}{2s+1}}\right)   &
\leq\frac{1}{\lambda}\left(  \frac{n}{\log n}\right)  ^{\frac{2s}{2s+1}%
}\mathbb{E}\left\Vert \widehat{\beta}_{\xi_{j,n}}-\beta\right\Vert _{2}^{2}\\
&  \leq\frac{M_{3}}{\lambda}\left(  \frac{n}{\log n}\right)  ^{\frac{2s}%
{2s+1}}\left(  \frac{n}{\log n}\right)  ^{-\frac{2s}{2s+1}}\leq M_{3}%
\lambda^{-1}.
\end{align*}

Hence, outside a set $\Omega_{n,2}$ of probability less than $M_{3}%
\lambda^{-1}$, it holds that $\left\Vert \widehat{\beta}_{\xi_{j,n}}%
-\beta\right\Vert _{2}^{2}\leq\lambda\left(  \frac{n}{\log n}\right)
^{-\frac{2s}{2s+1}}$. Therefore, the condition $\left\Vert \widehat{\beta
}_{\xi_{j,n}}-\beta\right\Vert _{2}\leq\frac{1}{2ebA_{j_{1}}}$ holds if
$\left(  \lambda\left(  \frac{n}{\log n}\right)  ^{-\frac{2s}{2s+1}}\right)
^{\frac{1}{2}}\leq\frac{1}{2ebA_{j_{1}}}$, which is equivalent to $4e^{2}%
b^{2}A_{j_{1}}^{2}\lambda\left(  \frac{n}{\log n}\right)  ^{-\frac{2s}{2s+1}%
}\leq1$. Using that $b^{2}\leq M_{1}^{2}e^{2\gamma_{j_{1}}+2\epsilon_{j_{1}}}$
the last inequality is true if $\eta_{j_{1},n}=4M_{1}^{2}e^{2\gamma_{j_{1}%
}+2\epsilon_{j_{1}}+2}A_{j_{1}}^{2}\left(  \frac{n}{\log n}\right)
^{-\frac{2s}{2s+1}}\leq\frac{1}{\lambda}$.

Hence, outside the set $\Omega_{n,2}$, our assumptions imply that there exists
$\widehat{\theta}_{n}=\theta\left(  \widehat{\beta}_{\xi_{j,n}}\right)  $ such
that $\left\langle f_{j_{1},\widehat{\theta}_{n},\xi_{j,n}}^{HT},\psi
_{j,k}\right\rangle =\widehat{\beta}_{\xi_{j,n},\left(  j,k\right)  }$ for all
$\left(  j,k\right)  \in\Lambda_{j_{1}}$. Finally, outside the set
$\Omega_{n,2}$, by using the bound given in Lemma \ref{Lema3.2enAB}, one
obtains the following inequality for the estimation error%
\[
\Delta\left(  f_{j_{1},\theta_{j_{1}}^{\ast}};f_{j_{1},\widehat{\theta}%
_{n},\xi_{j,n}}^{HT}\right)  \leq2M_{1}e^{\gamma_{j_{1}}+\epsilon_{j_{1}}%
+1}\lambda\left(  \frac{n}{\log n}\right)  ^{-\frac{2s}{2s+1}},
\]
which completes the proof of Proposition \ref{prop:existnonlin}.
\end{proof}

Our assumptions on $j_{1}(n)$ imply that $\frac{1}{2}\frac{n}{\log n}%
\leq2^{j_{1}(n)}\leq\frac{n}{\log n}$. Therefore, using Lemma \ref{LemmaAjDj},
one has that for all $f\in F_{p,q}^{s}(M)$ with $s>1/p$,
\[
\gamma_{j_{1}\left(  n\right)  }\leq C\left(  \frac{n}{\log n}\right)
^{-\left(  s-\frac{1}{p}\right)  },\quad A_{j_{1}\left(  n\right)  }%
\leq\left(  \frac{n}{\log n}\right)  ^{\frac{1}{2}},\quad D_{j_{1}\left(
n\right)  }\leq C\left(  \frac{n}{\log n}\right)  ^{-s^{\ast}},
\]
where $C$ denotes constants not depending on $g=\log(f)$. Hence,%
\[
\underset{n\rightarrow+\infty}{\lim}\epsilon_{j_{1}\left(  n\right)
}=\underset{n\rightarrow+\infty}{\lim}2M_{1}^{2}e^{2\gamma_{j_{1}\left(
n\right)  }+1}A_{j_{1}\left(  n\right)  }D_{j_{1}\left(  n\right)  }=0,
\]
uniformly over $F_{p,q}^{s}(M)$ for $s>1/p$. For all sufficiently large $n$,
$\epsilon_{j_{1}\left(  n\right)  }\leq1$ and thus, using Proposition
\ref{prop:existnonlin}, there exists $\theta_{j_{1}(n)}^{\ast}\in{\mathbb{R}%
}^{\#\Lambda_{j_{1}(n)}}$ such that
\[
\Delta\left(  f;f_{j_{1}(n),\theta_{j_{1}(n)}^{\ast}}\right)  \leq\frac{M_{1}%
}{2}e^{\gamma_{j_{1}(n)}}D_{j_{1}(n)}^{2}\leq C\left(  \frac{n}{\log
n}\right)  ^{-2s^{\ast}}\mbox{ for all }f\in F_{p,q}^{s}(M).
\]
Now remark that if $p=2$ then $s^{\ast}=s>1$ (by assumption), thus%
\[
\Delta\left(  f;f_{j_{1}(n),\theta_{j_{1}(n)}^{\ast}}\right)  =O\left(
\left(  \frac{n}{\log n}\right)  ^{-2s}\right)  \leq O\left(  \left(  \frac
{n}{\log n}\right)  ^{-\frac{2s}{2s+1}}\right)  .
\]
If $1\leq p<2$ then one can check that condition $s>\frac{1}{2}+\frac{1}{p}$
implies that $2s^{\ast}>\frac{2s}{2s+1}$, hence%
\begin{equation}
\Delta\left(  f;f_{j_{1}(n),\theta_{j_{1}(n)}^{\ast}}\right)  \leq O\left(
\left(  \frac{n}{\log n}\right)  ^{-\frac{2s}{2s+1}}\right)
.\label{deltanonlin1}%
\end{equation}

By the same arguments it holds that%
\[
\underset{n\rightarrow+\infty}{\lim}\eta_{j_{1}\left(  n\right)  ,n}%
=\underset{n\rightarrow+\infty}{\lim}4M_{1}^{2}e^{2\left(  \gamma
_{j_{1}\left(  n\right)  }+\epsilon_{j_{1}\left(  n\right)  }+1\right)
}A_{j_{1}\left(  n\right)  }^{2}\left(  \frac{n}{\log n}\right)  ^{-\frac
{2s}{2s+1}}=0,
\]
uniformly over $F_{p,q}^{s}(M)$ for $s>1/p$. Now let $\lambda>0$. The above
result shows that for sufficiently large $n$, $\lambda\leq\eta_{j_{1}\left(
n\right)  ,n}^{-1}$, and thus using Proposition \ref{prop:existnonlin} it
follows that there exists a set $\Omega_{n,2}$ of probability less than
$M_{3}\lambda^{-1}$ such that outside this set there exists $\widehat{\theta
}_{n}\in$ ${\mathbb{R}}^{\#\Lambda_{j_{1}(n)}}$ which satisfies:%
\begin{equation}
\Delta\left(  f_{j_{1}\left(  n\right)  ,\theta_{j_{1}\left(  n\right)
}^{\ast}};f_{j_{1}\left(  n\right)  ,\widehat{\theta}_{n},\xi_{j,n}}%
^{HT}\right)  \leq2M_{1}e^{\gamma_{j_{1}\left(  n\right)  }+\epsilon
_{j_{1}\left(  n\right)  }+1}\lambda\left(  \frac{n}{\log n}\right)
^{-\frac{2s}{2s+1}}\label{deltanonlin2}%
\end{equation}
for all $f\in F_{p,q}^{s}(M)$. Then, by the Pythagorian-like identity
(\ref{IdPyth}) it follows that outside the set $\Omega_{n,2}$,
\[
\Delta\left(  f;f_{j_{1}\left(  n\right)  ,\widehat{\theta}_{n},\xi_{j,n}%
}^{HT}\right)  =\Delta\left(  f;f_{j_{1}(n),\theta_{j_{1}(n)}^{\ast}}\right)
+\Delta\left(  f_{j_{1}\left(  n\right)  ,\theta_{j_{1}\left(  n\right)
}^{\ast}};f_{j_{1}\left(  n\right)  ,\widehat{\theta}_{n},\xi_{j,n}}%
^{HT}\right)  ,
\]
and thus Theorem \ref{Teo5.1enAB} follows from inequalities
(\ref{deltanonlin1}) and (\ref{deltanonlin2}).

\subsection{\textbf{Proof of Theorem \ref{Teo2enComte} }}

The proof is analogous to the one of Theorem \ref{Teo5.1enAB}. It follows from
Lemma \ref{LemmaBoundBiasVarThresBetahatEpsHat}.

\end{document}